\documentclass{article}
\usepackage{amsmath,amssymb,amscd,amsthm}
\usepackage[ansinew]{inputenc}
\usepackage[pdftex]{graphicx}
\usepackage{hyperref}
\usepackage{tikz,tikz-cd}

\newtheorem{thm}{Theorem}
\newtheorem{prop}{Proposition}

\newcommand{\rsp}{\raisebox{0em}[2.7ex][1.3ex]{\rule{0em}{2ex} }}
\newcommand{\OO}{\mathcal O}
\newcommand{\what}[1]{{#1}\,\widehat{\phantom{n}\!\!}} 
\newcommand{\Z}{{\mathbb Z}}
\newcommand{\Q}{{\mathbb Q}}
\newcommand{\R}{{\mathbb R}}
\newcommand{\cE}{{\mathcal E}}
\newcommand{\cN}{{\mathcal N}}
\newcommand{\cK}{{\mathcal K}}
\newcommand{\fm}{{\mathfrak m}}
\newcommand{\fp}{{\mathfrak p}}
\newcommand{\fP}{{\mathfrak P}}
\newcommand{\fq}{{\mathfrak q}}
\newcommand{\fM}{{\mathfrak M}}
\newcommand{\fN}{{\mathfrak N}}
\newcommand{\fg}{{\mathfrak g}}
\newcommand{\fa}{\mathfrak a}
\newcommand{\frf}{{\mathfrak f}}
\newcommand{\cH}{\mathcal H}
\newcommand{\ord}{\operatorname{ord}\,}
\newcommand{\Ann}{{\operatorname{Ann}}}
\newcommand{\Tak}{\operatorname{Tak}}
\newcommand{\Cl}{{\operatorname{Cl}}}
\newcommand{\GL}{{\operatorname{GL}}}  
\newcommand{\Gal}{{\operatorname{Gal}}}
\newcommand{\disc}{\operatorname{disc}}
\newcommand{\gen}{{\operatorname{gen}}}
\newcommand{\cen}{{\operatorname{cen}}}
\newcommand{\Hom}{{\operatorname{Hom}}}
\newcommand{\im}{{\operatorname{im}\,}}
\newcommand{\eps}{\varepsilon}
\newcommand{\lra}{\longrightarrow}
\newcommand{\rd}{\operatorname{rd}\,}
\newcommand{\la}{\langle}
\newcommand{\ra}{\rangle}
\newcommand{\Wt}[1]{\widetilde{#1}}

\title{Arnold Scholz: Between Mathematics and Politics}
\author{Franz Lemmermeyer}

\begin{document}

\maketitle

\begin{abstract}
  This is an English translation (with a few modifications) of the
  introduction to Scholz's life and work taken from the edition  \cite{LRS}
  of the correspondence between Helmut Hasse, Arnold Scholz and Olga Taussky.
  I hope to find the time one day to include more information on
  Scholz's mathematical results.
\end{abstract}

In this article we want to present, apart from a brief
biography,\footnote{N.~Korrodi has kindly helped me to fill some gaps.}
an introduction to Arnold Scholz's number-theoretic work and
at the same time to convey an impression of the difficulties he faced
from the National Socialists' seizure of power until his death.

Arnold Scholz was born on 24 December 1904 in the Charlottenburg
district of Berlin.  His father, Reinhold Scholz (1857--1933), was
employed as a physicist and mathematician at the Military Research
Institute in Berlin; his mother was Johanna Scholz, née Diesfeld
(1874--1957).

In 1923 Scholz passed the Abitur at the Kaiserin-Augusta-Gymnasium in
Charlottenburg and then attended lectures in philosophy, mathematics
and musicology at the University of Berlin.  Scholz spent the summer
of 1927 abroad in  Vienna studying with Philipp Furtwängler; he received
his Ph.D. in 1928 \emph{magna cum laude} under Issai Schur and moved
to Freiburg in April 1929.

Arnold Scholz quickly became one of the foremost experts on Takagi's
class field theory; what others had to prove, he simply ``saw''. This
is one reason why his publications are sometimes very difficult to read.
According to Taussky's obituary of Scholz, Schur is said to have
remarked about one of Scholz's early papers: ``More Landau, less
Goethe!''

A quite different difficulty is described by Scholz himself in his letter
from 19 December 1935:
\begin{quote}
  \em The final formulation is so difficult for me with a new subject
  because I usually do not think in words, but in spatial images,
  furnished with colours and all sorts of other associations.
\end{quote}

Hasse, in an expert opinion on Scholz, criticised him for using too
little modern algebra:
\begin{quote}
  \em On the other hand, this strongly developed sense for numbers
  is also a certain weakness, for it often prevents him from fully
  employing the meanwhile developed more formal methods of modern
  algebra, which in recent years have so decisively influenced
  number theory.
\end{quote}

In some sense this may well be true, even though the readability
of Scholz's articles suffers primarily from the fact that he does not
condense recurring patterns of arguments into lemmas and propositions,
but rather hides them in the text at their first appearance and then
employs them a second time without any comment: too much Goethe, too
little Landau.

Concerning modern algebra, the following remark by Jehne from
\cite[S.~221]{Jehne} seems to me more to the point:
\begin{quote}
  \em Without any cohomology theory at hand Scholz recognized that
  the number knot $\nu_{K|k}$ is an epimorphic image of the fundamental
  group $\pi \fg$ of the Galois group in a natural way.
\end{quote}

It is quite simply the case that the language of exact sequences and
cohomology theory -- so helpful in questions of embedding problems and
norm residues -- had not yet been developed, and Scholz was therefore
unable to express his results with the desired clarity using the
language available to him at the time.

This also explains why most of his results had to be rediscovered
first before one could find them in Scholz's papers. One of the few
who studied Scholz's works thoroughly was Shafarevich: by extending
Scholz's methods for constructing number fields with given Galois
group, he was able to show that every solvable finite group occurs as
a Galois group over $\Q$ (Scholz had proved this for $p$-groups and
for certain two-step groups, i.e.\ those whose commutator subgroup is
abelian).

Other results of Scholz that were rediscovered much later include:
\begin{itemize}
\item The theory of genus and central class fields, initiated by
  Chebotarev and extended by Scholz, was developed a second time
  by Fröhlich \cite{FrCE}.
\item The ``Scholz reciprocity law'', already proved by Schönemann,
  was rediscovered by Emma Lehmer (see \cite{Lehmer}) and only then
  located in Scholz's work; it triggered a flood of investigations
  into power residue characters of units. By contrast, Scholz's
  generalization of this law to primes $\ell > 2$ has so far played no
  role in the literature.
\item Scholz's theory of knots was modernized by Jehne \cite{Jehne}
  (after Tate had rediscovered the main theorem in cohomological
  form) and subsequently further developed by many authors.
\item Parts of the work of Scholz \& Taussky \cite{S17} were presented
  in modern form by Heider \& Schmithals \cite{HS82}.
\item The theory of abelian crossings, developed by Scholz in
  connection with his knots, was further elaborated by Heider
  \cite{HeiderK} (cf.\ also Steinke \cite[\S~2]{Steinke}) using modern
  tools.
\end{itemize}

Similar investigations are still pending for other works of Scholz;
even as regards the theory of knots, Jehne's paper represents only a
first step towards a proper ``translation'' of Scholz's work: since
Jehne we know what Scholz proved, but his proofs have remained
obscure to this day.

In the following we wish to go through the most important stations in
the life and work of Arnold Scholz in chronological order.

\section*{1927}
Together with his former teacher Friedrich Neiß and his fellow student
Klaus Müller, Scholz studied Hasse's class field report and, in his
first letter to Hasse dated 22 April 1927, made suggestions for
simplifying several proofs, which ultimately led to their joint paper
\cite{S1}. Scholz's essential contributions are the
following: He realised that, starting from Takagi's definition of
class fields, one can show directly and quickly
\begin{itemize}
\item that class fields are necessarily Galois and possess a Galois
  group all of whose subgroups are normal;
\item that class fields possess simple functorial properties: if
  $L/K/k$ is a tower of field extensions
  and $L/k$ is a class field, then so are $L/K$ and $K/k$.
\end{itemize}

The proofs are of an analytic nature and make essential use of the
pole of the Dedekind zeta function at $s = 1$. With a little group
theory, one easily obtains from these statements a proof that class
fields are {\em abelian} extensions. Some of Scholz's ideas from 
1927 were only published later in the paper \cite{S10}.

On 23 April Scholz travelled to Vienna and remained there until
September in order to study with Furtwängler the topic of his
dissertation, namely the construction of field extensions
with prescribed solvable Galois group using class field theory. Later
others also began to take an interest in this topic, in particular
Richter \cite{Rina,Riab}, Reichardt \cite{ReiKZk}, and Tannaka
\cite{Tannaka}. Definitive results were obtained much later by
Shafarevich \cite{Shaf}, who built upon Scholz's ideas.

In order to explain the basic idea of Scholz, we shall carry out his
construction for a few very simple examples.

\subsection*{Octic dihedral extensions}
Let
$$ G = D_4 = \langle S, T : S^4 = T^2 = 1,\ TST = S^{-1} \rangle $$ be
the dihedral group of order $8$. The construction of extensions $K/\Q$
with $\Gal(K/\Q) \simeq D_4$ is of course no problem that would
require class field theory for its solution; our aim here is to
illustrate what the Scholz construction looks like in such a simple
example.

Scholz constructs such two-step extensions in two abelian steps; in
the present case there are three possibilities:
\begin{enumerate}
\item[a)] construct $K$ as a suitable quadratic
          extension $K/k$ of the biquadratic extension
          $k = \Q(\sqrt{m},\sqrt{n}\,)$;
\item[b)] construct $K$ as a suitable biquadratic extension
          of a quadratic subfield;
\item[c)] construct $K$ as a cyclic extension $K/k$ of degree four
          over a quadratic subfield $k = \Q(\sqrt{m}\,)$.
\end{enumerate}

We shall first choose the third variant here; in order to be able to
describe the Galois group of the extension $K/\Q$ precisely at each
step, we wish to realise the generating automorphisms in each step as
inertia automorphisms of suitable prime ideals, and for this purpose
it is helpful if exactly one prime ideal ramifies in each step. In the
present case we choose a quadratic extension $k/\Q$ in which exactly
one prime $p$ ramifies. By the Kronecker--Weber theorem these fields
are the quadratic subfields of the field of $p$-th roots of unity; in
order that the quadratic subfield does not ramify at the infinite
places (which would introduce an additional difficulty for continuing
the construction), we must choose $p \equiv 1 \pmod{4}$
(alternatively: the unit $-1$ in the base field should be a quadratic
residue modulo $p$; we shall encounter such conditions quite
frequently). Thus $k = \Q(\sqrt{p}\,)$ is a real quadratic subfield
whose Galois group is generated by the Frobenius automorphism $T$ of
the prime ideal $(\sqrt{p}\,)$.

\begin{figure}[ht!]
  \begin{minipage}{6cm}
    \begin{tikzpicture}      
  \node (Q)  at (0,0)   {$\Q$};
  \node (k)  at (0,1.5)   {$k$};
  \node (K1) at (-2,3)  {$K_1$};
  \node (K) at (-0.8,3)  {$K$};
  \node (Kt) at (0.8,3)  {$\Wt{K}$};
  \node (K2) at ( 2,3)  {$K_2$};
  \node (L)  at (0,4.5)  {$L$};
  \draw (Q) -- (k) -- (K1) -- (L) -- (K2) -- (k) -- (K) -- (L) -- (Kt) -- (k);
\end{tikzpicture} \end{minipage} \quad \begin{minipage}{6cm}
\begin{tikzpicture}
  \node (Q)  at (0,0)    {$\Gamma$};
  \node (k)  at (0,1.5)  {$\la S_1, S_2 \ra$};
  \node (K1) at (-2,3)   {$S_1$};
  \node (K) at (-0.8,3)  {$S_1S_2$};
  \node (Kt) at (0.8,3)  {$S_1S_2^{-1}$};
  \node (K2) at ( 2,3)   {$S_2$};
  \node (L)  at (0,4.5)  {$1$};
  \draw (Q) -- (k) -- (K1) -- (L) -- (K2) -- (k) -- (K) -- (L) -- (Kt) -- (k);
\end{tikzpicture}
\end{minipage}
\caption{Extension with wreath product $C_4\,\wr\,C_2$}\label{HasKP}
\end{figure}
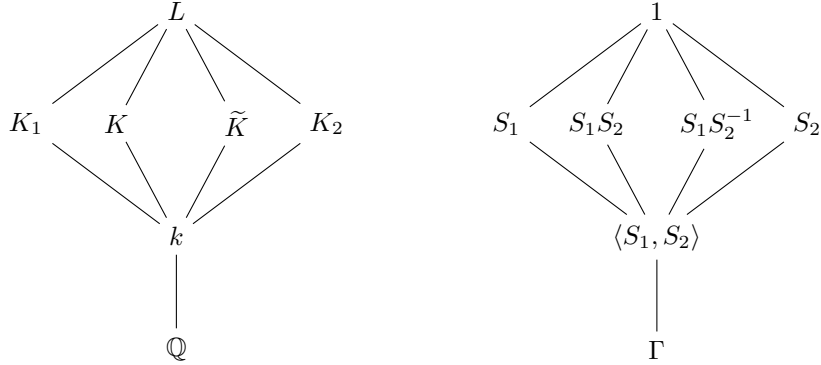

Now we turn to the construction of the cyclic extension $K/k$. In
order that $T$ retains order $2$ also in $K/\Q$, we shall choose $K$
in such a way that the prime $p$ does not ramify further in $K/k$;
then the inertia group of $p$ in $K/\Q$ will still have order $2$ and
be generated by an element $\widetilde{T}$ of order $2$ whose
restriction to $k$ is $T$. In what follows we shall, like Scholz, not
distinguish between $\widetilde{T}$ and $T$.  The
construction of a cyclic extension $K/k$ with the aid of class field
theory could certainly also be carried out directly; in any case, it
is not sufficient to lift cyclic extensions from $\Q$ to $k$, since
the Galois group of such an extension would necessarily be abelian. In
order to ensure that the cyclic extension $K/k$ does not come from a
cyclic extension $F/\Q$, Scholz proceeds via the construction of an
extension $L/\Q$ whose Galois group is the wreath product (at Scholz
this construction is called, as we shall see below, a ``disposition
group'')
\begin{align}
  \notag \Gamma & = C_4\,\wr\,C_2
       = \langle S_1, S_2, T : S_1^4 = S_2^4 = T^2 = 1, \\
\label{Cwr24} & \qquad S_1S_2 = S_2S_1,\ TS_1T = S_2,\ TS_2T = S_1 \rangle
\end{align}
here and below we denote by $C_m$ the cyclic group of order
$m$. In this case, the construction is also carried out in such a way
that exactly one prime ideal ramifies in each of the cyclic base
extensions. If one introduces in $\Gamma$ the additional relation
$S_1 = S_2$ (or $S_1 = S_2^{-1}$), i.e.\ forms the quotient group
$\Gamma/\langle S_1S_2^{-1} \rangle$ (or $\Gamma/\langle S_1S_2 \rangle$),
then the fixed field $K$ of the normal subgroup $\langle S_1S_2^{-1} \rangle$
(or $\langle S_1S_2 \rangle$) is a cyclic
extension of $k$ in which two prime ideals ramify; in the first case
we have $\Gal(K/\Q) \simeq C_2 \times C_4$, in the second $\Gal(K/\Q)
\simeq D_4$ (the Hasse diagram in Fig.~\ref{HasKP} does not contain all
subfields or subgroups).

Scholz therefore wishes to construct a cyclic extension $K_1/k$ that
is conjugate over $k$ to a cyclic extension $K_2/k$ with $K_1 \cap K_2 = k$.
This situation arises automatically if $K_1/k$ is a cyclic
extension whose conductor $\frf_1 = \fq$ is a prime ideal of degree 1;
by class field theory, the corresponding class field $K_2$ to the
conductor $\frf_2 = \fq^T$ will then also be cyclic of degree 4 and
disjoint from $K_1$ over $k$.

To construct the class field $K_1/k$, we therefore choose a prime
number $q$ that splits in $k/\Q$: $q\OO_k = \fq_1\fq_2$. The ray class
number for an extension with conductor $\frf_1 = \fq_1$ is
$$ h_k(\fq_1) = h_k \cdot \frac{\Phi(\fq_1)}{(E:E^{(1)})}, $$
where $h_k$ is the class number of $k$, $\Phi(\fq_1) = q-1$ is the
number of residue classes in $k$ coprime to $\fq_1$, $E = \langle -1,
\eps_p \rangle$ is the unit group of $k$, and $E^{(1)}$ denotes the
subgroup of all units $\eps \equiv 1 \pmod{\fq_1}$.

\medskip\noindent
The minimum requirement on $q$ is therefore that this ray class number
is divisible by $4$. For this to hold, $q-1$ must be a multiple
of $4$, and in order for $E^{(1)}$ to have index divisible by $4$ in
$E$, both $-1$ and the fundamental unit $\eps_p$ must be quadratic
residues modulo $\fq_1$. If we choose $q \equiv 1 \pmod{4}$, then
$q-1$ is a multiple of $4$ and $-1$ is a quadratic residue modulo
$\fq_1$; the question of when $\eps_p$ is a quadratic residue modulo
$\fq_1$ was later studied extensively by Scholz in \cite{S15}.

The existence of such primes $q$ can easily be proved: we need only
choose $q$ so that this prime splits completely in the extension
$\Q(\sqrt{p}, \sqrt{-1}, \sqrt{\eps_p}\,)$, which is always possible
by classical density theorems (Kronecker, Frobenius, later
Chebotarev).

By this choice of $q$ we have now obtained a cyclic extension
$K_1/k$ of degree $4$ in which at most the prime ideal above $\fq_1$
ramifies. In the present case this prime ideal {\em must} ramify,
because $k$ possesses no quadratic unramified extensions. Let $K_1$
therefore be the cyclic subextension of degree $4$ in the ray class
field of $k$ modulo $\fq_1$, and let $S_1$ be the inertia automorphism
of $\fq_1$; then $S_1^4 = 1$, and by Artin the extension $K_2/k$
conjugate over $k$ to $K_1/k$ is also cyclic, its Galois group is
generated by the inertia automorphism $S_2$ of $\fq_2$, and the
compositum $L = K_1K_2$ has Galois group as in (\ref{Cwr24}).

The desired dihedral extension is obtained by setting $S_1 = S_2^{-1}$
in (\ref{Cwr24}), i.e.\ by considering the quotient group of $\Gamma$
by $\langle S_1S_2 \rangle$: this normal subgroup then fixes a
dihedral extension $K/\Q$ that is cyclic over the quadratic subfield
$k = \Q(\sqrt{p}\,)$; in $K/k$ precisely the two prime ideals above
$q$ will ramify.

We now also see that constructing a dihedral extension by the second
method would have been much simpler, since the dihedral group $D_4$ is
isomorphic to the wreath product $C_2\,\wr\,C_2$, and we would
therefore only need to construct an extension of a quadratic number
field as a subfield of the ray class field modulo a primary prime
ideal.

\subsection*{The ramification problem}
In our above construction of a dihedral extension we were, at one
point, somewhat fortunate in that the base field $k = \Q(\sqrt{p}\,)$
has odd class number. In general we cannot assume that our base field
satisfies such arithmetic conditions. This fact leads to problems that
we wish to illustrate again with a very simple example.

Let $k$ be an imaginary quadratic number field different from
$\Q(\sqrt{-1}\,)$ and $\Q(\sqrt{-3}\,)$, and let $p > 2$ be a prime
that splits in $k$: $p\OO_k = \fp_1\fp_2$. The ray class number of $k$
modulo $\fp$ is then
$$ h_k\{\fp\} = h_k \cdot \frac{p-1}{2}. $$
The class number $h_k$ is even if and only if there exist unramified
quadratic extensions of $k$. For the existence of a quadratic
extension ramified exactly at $\fp$, however, it is necessary that $p
\equiv 1 \pmod{4}$. This condition is usually also sufficient: for the
prime ideal $\fp = (3 + 2 \sqrt{-5}\,)$ of norm $29$ in the quadratic
number field $k = \Q(\sqrt{-5}\,)$ with class number $h_k = 2$, we
have $h_k\{\fp\} = 28 \equiv 0 \pmod{4}$, but the degree-$4$ subfield
of the ray class field modulo $\fp$ consists of the unramified
extension $K = k(i)$, the Hilbert class field, and a quadratic
extension $L = K\big(\sqrt{(1+2i)(3+2\sqrt{-5}\,)}\big)$ adjoined to
it, in which $\fp$ ramifies.

There is, however, no quadratic extension of $k$ that is ramified
exactly at $\fp$; the ``usual suspects'' such as
$k(\sqrt{\pm 3 + 2 \sqrt{-5}}\,)$ are in fact also ramified over
the prime ideal above $2$.

If, therefore, in a given base field $k$ one has found a prime ideal $\fp$
for which the ray class number $h_k(\fp)$ is divisible by $\ell$, then there
exists a cyclic extension $K/k$ of degree $\ell$ in which at most $\fp$
ramifies; however, this extension may also be unramified, which leads to
insurmountable problems when continuing the construction, because one can
say little about the corresponding automorphisms, as they are tied to the
behaviour of the associated ideal classes in extensions. If one chooses $\fp$
in such a way that $h_k(\fp)/h_k$ is divisible by $\ell$, then there does
indeed exist an extension $F/H$ of the Hilbert class field $H$ of $k$ that
is ramified exactly over $\fp$, but in general this extension cannot be
pulled back down, except when the ray class group is the {\em direct product}
of the ideal class group and a group of order divisible by $\ell$.

Together with Furtwängler, Scholz discovers how one can achieve this
decomposition of the ray class group by a clever choice of the prime ideal
$\fp$. To this end, let $K$ be a number field containing a primitive
$\ell$-th root of unity $\zeta$ (what follows also holds more
generally for roots of unity of prime power degree); a prime ideal
$\fp$ is called $\ell$-primary if it is coprime to $\ell$ and if
$(\omega/\fp)_\ell = 1$ for all singular $\omega \in K \setminus
(1-\zeta)$, i.e., all elements $\omega$ coprime to $\ell$ that are
$\ell$-th ideal powers: $(\omega) = \fa^\ell$. In fields with class
number coprime to $\ell$, $\fp$ is $\ell$-primary if and only if every
unit is an $\ell$-th power residue modulo $\fp$. The Frobenius density
theorem guarantees the existence of such prime ideals.

Scholz observes that he can also use the notion of primary prime ideal
in number fields that do not contain the corresponding roots of unity:
he calls a prime ideal $\fp$ in $K$ {\em primary} modulo $\ell^h$ if
$N\fp \equiv 1 \pmod{\ell^h}$ and all elements of $K$ that are
$\ell^h$-th ideal powers are also $\ell^h$-th power residues modulo
$\fp$.

In the case $k = \Q(\sqrt{-5}\,)$, the prime ideal $\fp$ above $29$ is
not $2$-primary, because the singular number $\omega = 2 + \sqrt{-5}$
is a quadratic non-residue modulo $\fp$. By contrast, the prime ideal
$\fq = (-3+4\sqrt{-5}\,)$ of norm $89$ is $2$-primary; since $N\fq
\equiv 1 \pmod{4}$, we have $(\frac{-1}{\fq}) = +1$, while $4
\sqrt{-5} \equiv 3 \pmod{\fq}$ implies $2 + \sqrt{-5} \equiv
\frac{1}{4}(8 + 4 \sqrt{-5}) \equiv \frac{11}{4} \equiv 25 \pmod{\fq}$
and thus certainly $(\frac{2+\sqrt{-5}}{\fq}) = +1$. Indeed,
$k(\sqrt{-3+4\sqrt{-5}}\,)$ is a quadratic extension in which
precisely the prime ideal $\fq$ ramifies.

For a prime ideal $\fp$ that is primary modulo $\ell^h$, consider 
the ideal group $\cH$ consisting of all ideals $\fa$ in $k$ coprime to
$\fp$ such that some $m$-th power with $\ell \nmid m$ equals a
principal ideal generated by an $\ell^h$-th power residue $\alpha$:
$\fa^m = (\alpha)$ with $(\alpha/\fp)_{\ell^h} = 1$. This last
condition does not depend on the choice of $\alpha$, because every
unit is an $\ell^h$-th power residue modulo $\fp$. If $D_k$ denotes
the group of all ideals coprime to $\fp$, then
$\Cl_\ell\{\fp\} = D_k/\cH$ is a finite $\ell$-group that has the
$\ell$-class group as a quotient group. More precisely, Scholz shows that
when $\fp$ is chosen to be primary, we always have
$\Cl_\ell\{\fp\} \simeq \Cl_\ell(k) \times A$;
the class field corresponding to the quotient group
$D_k/\Cl_\ell(k)$ is what Scholz calls a {\em residue field}, i.e., a
purely ramified cyclic extension of degree $\ell^h$ in which exactly
the prime ideal $\fp$ (and only this one) is purely ramified.

\begin{figure}[p!]
\noindent
\includegraphics[width=11.5cm]{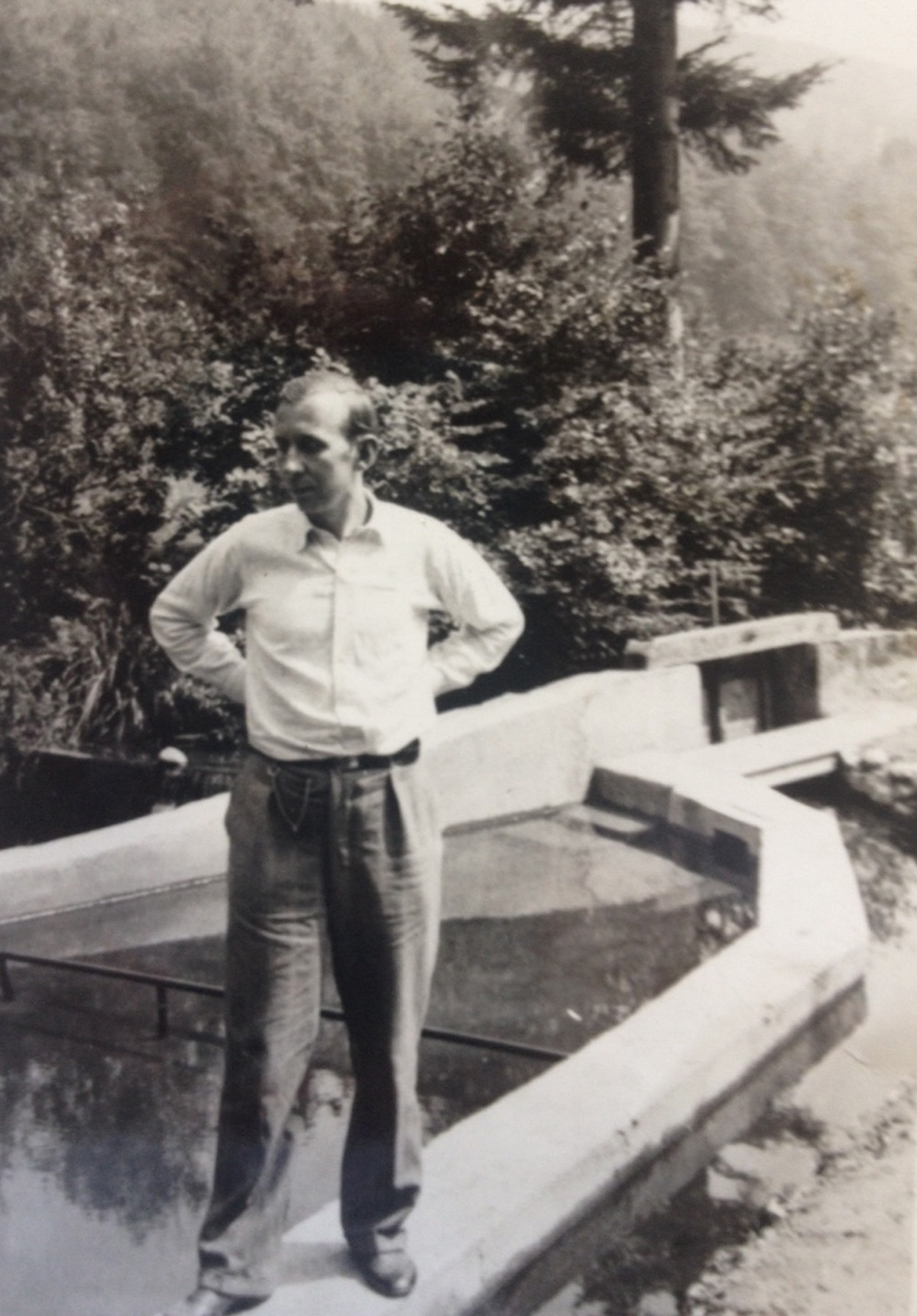}
\caption{Arnold Scholz}\label{Sch1}
\end{figure}

\subsection*{Disposition groups}
If $G$ is a solvable group, then extensions with Galois group $G$ can
be constructed step by step in abelian steps. Scholz recognised that
it suffices, for the construction of these abelian extensions, to
consider two types of groups: abelian $p$-groups and ``disposition
groups''. In \cite{S2} Scholz treats ``two-step disposition groups''
$G$; these are groups possessing the following two properties:
\begin{enumerate}
\item $G$ contains an abelian subgroup $A$ such that $G/A$ is abelian;
\item $A$ contains an element $a$ such that the elements $a^\sigma$
  form a basis of $A$ as $\sigma$ runs through a system of
  representatives of $G/A$.
\end{enumerate}
In modern terminology, a two-step disposition group is therefore
simply the wreath product $C_m\,\wr\,N$, where $C_m \simeq \langle a
\rangle$ and $N = G/A$. The simplest examples of wreath products are
the dihedral group $D_4 \simeq C_2\,\wr\,C_2$ and the alternating
group $A_4 = C_2\,\wr\,C_3$. For two-step disposition groups Scholz
can realise the step-by-step construction with exactly one ramified
prime ideal in each step.

\subsection*{Non-abelian extensions of degree $pq$}
Let $p$ and $q \equiv 1 \pmod{p}$ be two primes. Then there exists a
non-abelian group of order $pq$, namely
$$ \Gamma = \langle S, T : S^p = T^q = 1,\ S^{-1}TS = T^r \rangle , $$
where $r^p \equiv 1 \pmod{q}$ and $r \not\equiv 1 \pmod{q}$.

Following our preliminary considerations in the case of the dihedral
group, the procedure should now be clear: we construct a cyclic
extension $k/\Q$ of prime degree $p$ as a subextension of the field of
$\ell$-th roots of unity, where $\ell \equiv 1 \pmod{p}$ is prime. Then
we choose a $q$-primary prime ideal $\fq$ in $k$ and construct a
cyclic extension of prime degree $q$ with conductor $\fq$; the normal
closure of this extension then has Galois group $C_q\,\wr\,C_p$. The
desired extension with Galois group $\Gamma$ is finally obtained by a
clever formation of quotient groups.

Scholz shows in \cite{S2} how one can also construct the extension
with Galois group $\Gamma$ without the detour via the wreath
product. To this end, he introduces (in \cite[p.~353, bottom]{S2})
the notion of $(S-1)$-th power residues,
which plays a quite central role in a large number of his papers but
has received hardly any attention to this day.

To do so, Scholz considers a number field $K$ containing a primitive
$\ell$-th root of unity $\zeta$ and calls $\alpha \in \OO_K$ an
$(\ell,S-a)$-th power residue modulo a prime ideal $\fp$ (here we
tacitly assume that $\fp$ does not divide $\ell$ and that $\alpha$ is
coprime to $\fp \ell$) if
$$ \alpha \equiv \beta^{S-a} \gamma^\ell \pmod{\fp} $$
holds. More generally, for any element $\Lambda \in \Z[G]$, where
$\Z[G]$ denotes the group ring of the Galois group of $K/\Q$, one can
call an $\alpha \in \OO_K$ an $(\ell,\Lambda)$-th power residue modulo
$\fp$ if
$$ \alpha \equiv \beta^{\Lambda} \gamma^\ell \pmod{\fp} $$ holds. In
his papers Scholz applies this notion mainly to $\Lambda = (S-1)^2$.
In the present case, the aim is to adjoin to a cyclic extension of
degree $p$ an extension $L/K$ with Galois group $\Gamma$. To this end,
Scholz observes that the symbolic order $\fM = \Ann(A)$, i.e.\ the
annihilator of $A$ in the group ring $\Z[G]$, is in the present case
equal to $\fM = (q, S-r)$ and that this ideal in turn divides $(q,
S^p-1)$.

As before, Scholz chooses a $q$-primary prime ideal $\fq$, but he does
not construct an extension whose residue field is the full disposition
group $C_q\,\wr\,C_p$ corresponding to an ideal group whose principal
ideals consist of the $q$-th power residues modulo $\fq$; instead he
constructs the extension directly by including all $(q,S-a)$-th powers
in the ideal group.

Over the years a group-theoretic language has emerged that allows us
to state Scholz's main result from the paper \cite{S2} briefly and
concisely as follows: given a field extension $K/k$ with Galois group
$G$, and given the embedding problem in which one seeks an extension
$L/K$ such that $\Gal(L/k) = \Gamma$ and $\Gal(L/K) = N$, then this
embedding problem is always solvable whenever the exact sequence
$$ \begin{CD} 1 @>>> N @>>> \Gamma @>>> G @>>> 1 \end{CD} $$
splits.

Of the many smaller theorems scattered throughout the text of
\cite{S2}, we mention only one: if $K/k$ is normal with Galois group
$G$ and $L/K$ is a class field corresponding to the ideal group $H$,
then $K/k$ is normal if and only if $H$ remains fixed under the action
of all $\sigma \in G$.

\bigskip

After his return to Berlin, Scholz turns his attention to cubic number
fields with given discriminant and the related search for quadratic
number fields whose $3$-class group has rank $\ge 2$. These investigations
were later subsumed into the joint paper with Olga Taussky on the
capitulation problem. 

\section*{1928}
In February Scholz completes his dissertation on the construction of
extensions $K/K_0$ of a given number field $K_0$ whose Galois group is
a ``disposition group''. He now turns to the group-theoretic question
for which other groups one must solve the inverse Galois problem in
order to bring all solvable finite groups under control.

On 1 April Scholz takes up his assistant position with Erhard Schmidt
in Berlin. In August Scholz begins investigating the unit group and
ideal class group of number fields $\Q(\sqrt{p}, \sqrt{q}\,)$ for
primes $p \equiv q \equiv 1 \pmod{4}$, which will ultimately lead him
to questions concerning the solvability of the non-Pell equation $t^2
- pqu^2 = -4$, to the ``Scholz reciprocity law'', and also to his
counterexample to the conjecture that the local-global principle for
norms holds in non-cyclic extensions of number fields.

The corresponding investigations of composita of cyclic fields of
prime degree $\ell$ lead him to the discovery that there exist class
field towers of arbitrarily great length.

\section*{1929}
In 1929, Scholz receives his PhD in Berlin under Schur; as his topic,
he had chosen the construction of extensions of $\Q$ with prescribed
Galois group (see \cite{S2}). The ``inverse Galois problem'' remains a
central theme for Scholz, and his articles move back and forth between
class field theory on the one hand and group theory on the other e.g. 
in \cite{S3,S4}.

In \cite{S3} Scholz shows that the construction of number fields with
arbitrary two-step Galois group can be reduced to the construction of
extensions whose Galois group is a disposition group or a branch
group. Branch groups are here a family of, in a certain sense, maximal
$p$-groups that contain all $p$-groups as quotient groups.

\subsection*{Branch groups}
A two-step branch group of type $(\ell^{h_1},\ldots, \ell^{h_m};\ell^k)$
is defined by Scholz as a metabelian $\ell$-group $G$ with commutator
subgroup $A$ and the following properties:
\begin{enumerate}
\item $G/A$ has a basis $S_1 A$, \ldots, $S_m A$.
\item $\ord(S_\mu) = \ord(S_\mu A) = \ell^{h_\mu}$ for $\mu = 1,
  \ldots, m$. In particular, the cyclic groups generated by the
  $S_\mu$ intersect $A$ only at the identity element.
\item Between the commutators $a_{\mu\nu} = [S_\nu,S_\mu]$ with $\mu < \nu$
  there exist only the relations that follow from their definition:
  $$ a_{\lambda \mu}^{S_\kappa-1}
     a_{\mu \kappa}^{S_\lambda-1}
     a_{\kappa \lambda}^{S_\mu-1} = 1 $$     
     and                  
     $$ a_{\kappa \lambda}^{f_\kappa} = a_{\kappa \lambda}^{f_\lambda}= 1, $$
     where $f_\nu = \sum_{\lambda = 0}^{\ell^{h_\nu}-1} S_\nu^\lambda$.
     In particular, $a_{\mu\nu}$ has symbolic order
     $$ \fN_{\mu\nu}^{(k)} = \Big(\ell^k, \frac{S_\mu^{\ell^{h_\mu}}-1}{S_\mu-1}, 
          \frac{S_\nu^{\ell^{h_\nu}}-1}{S_\nu-1}, \ldots,
            S_\rho^{\ell^{h_\rho}}-1, \ldots (\rho \ne \mu,\nu)\Big). $$
\end{enumerate}
In \cite{S3} Scholz proved that every two-step metabelian $\ell$-group
with $m$ generators and commutators $a_{\mu\nu}$ whose orders divide
$\ell^k$ is a quotient group of such a branch group.

The simplest branch group has type $(2,2;2^k)$. Here $G$ is generated
by two elements $S_1$ and $S_2$ of order $2$ and is therefore a
dihedral group; the commutator $a = [S_1,S_2]$ has order $2^k$ and
generates the cyclic group $A = G'$, and we have $G/A \simeq C_2 \times C_2$.

\subsection*{Scholz in Freiburg}
After receiving his doctorate, Scholz moved to Freiburg in April 1929
(and not only in 1930, as Taussky writes in her obituary; see Remmert
\cite{RemmFB}); there he succeeded Baer (who had gone to Hasse in
Halle) as Loewy's assistant and met Zermelo, who was interested in the
foundations of mathematics and set theory.

Alfred Loewy was born on 20 June 1873 in Rawitsch near Posen and came
from a strictly orthodox Jewish family. He studied in Breslau, Munich,
Berlin and Göttingen and received his doctorate in 1894 under
Lindemann in Munich. In 1897 he received his habilitation in Freiburg,
became extraordinary professor in 1902 and full professor in 1919,
succeeding Stickelberger.

In 1916 Loewy went blind in one eye; a failed operation in 1928 cost
him the second eye as well. His doctoral student Krull (PhD 1921)
introduced modern algebra in the sense of Emmy Noether in
Freiburg. Another of Loewy's doctoral students, namely F.K. Schmidt,
received, after his doctorate, an assistant position created at the
request of Heffter and Loewy; as his successor Loewy brought Reinhold
Baer to Freiburg in 1926, and after the latter's habilitation the
position went to Arnold Scholz.

In 1933, Loewy is forced to retire by the National Socialists. Because
the decree of April 1933 is temporarily suspended, Loewy is still able
to give a lecture in the summer semester of 1933. Efforts by Heffter
to obtain an exceptional permit for Loewy fail, just as similar
initiatives elsewhere fail for, among others, Emmy Noether, Issai
Schur, and Edmund Landau.

Loewy's succession is arranged by Doetsch; since algebra and number
theory are said to have ``taken a completely abstract direction'', and
indeed ``under strong influence from the Jewish side'', a geometer is
to be appointed. The choice falls on Wilhelm Süss, who, with his
diplomatic skill, knows how to exploit the National Socialist moment
for his own purposes like few others.

With Scholz's departure, algebra is dead in Freiburg. Loewy writes his
book \emph{Blinde große Männer} (``Blind Great Men''), for which he
finds no publisher in National Socialist Germany: the book appears
posthumously in Zurich. Loewy dies on 25 January 1935 following an
operation in Freiburg. As Heffter writes, at his funeral, apart from a
few of his mathematical colleagues, only he himself and Hans
Spemann\footnote{Hans Spemann (1869--1941) was a biologist, active in
  Freiburg from 1919 until his retirement in 1937, and rector from
  1923 to 1924. He received the Nobel Prize in Physiology or Medicine
  in 1935.} attend. His wife followed him a few years later by
suicide.

Ernst Zermelo was born on 27 July 1871 in Berlin.  After studies in
Berlin, Halle, and Freiburg, he received his doctorate in Berlin in
1894 with a thesis on the calculus of variations. In 1897 he
habilitated in Göttingen and achieved fame in 1904 through his
clarification of the axiom of choice, which earned him a professorship
in Göttingen in 1905.

In 1910 he accepted a professorship in Zurich, which he gave up in
1916 for health reasons and moved to the Black Forest.

In 1926, on the application of Heffter and Loewy, he received an
honorary professorship in Freiburg, which he held until his
denunciation in 1935 by Doetsch's assistant Eugen Schlotter (the full
text of this denunciation can be found in \cite[pp.~296--297]{Ebb}):
Zermelo refused to begin his lectures with the Hitler salute.

During their common time in Freiburg, Arnold Scholz became Zermelo's
closest friend until Scholz's death in 1942, which deeply affected
Zermelo. Zermelo died on 21 May 1953 in Freiburg.
For Zermelo's life and work, see \cite{Ebb}.

\bigskip\noindent

In \cite{S5}, Scholz proves that there exist class field towers of
arbitrarily great length. At the end of the 19th century, Hilbert had
still conjectured in \cite{Hil2a,Hil2,Hil3} that the class number of
the Hilbert class field of number fields with class number $4$ is
always odd. This was disproved by Furtw\"angler in 1916 in \cite{Fw16}
using counterexamples; since then, however, no further progress had
been made on this question.

Scholz now succeeded in proving that a cyclic number field of prime
degree $\ell$ possesses an $\ell$-class field tower of any prescribed
minimal length, provided its $\ell$-class group has sufficiently large
rank.

The definitive answer to the question of infinitely long class field
towers was provided only by Golod and Shafarevich; refinements of the
question remain the subject of numerous investigations even today
(Benjamin, Lemmermeyer \& Snyder \cite{BLS1,BLS2,BLS3,BLS4}; Azizi et
al.\ (see e.g.\ \cite{Azizi}, \cite{Mouhib}); Bartholdi \cite{BaBu};
Bush \cite{Bush2,Bush,BM}; Hajir \cite{Hajir1,Hajir2}; Kuhnt
\cite{Kuhnt}; Maire \cite{Maire,MaML}; D.~Mayer
\cite{May12a,May12b,May13,May14}; McLeman \cite{McL}; Nover
\cite{Nover}; Steurer \cite{Steurer}).
The group theoretic side of this topic has also been
thoroughly investigated (see Magnus \cite{Magnus}, Serre \cite{Ser},
Nebelung \cite{Neb89}, and more recent articles e.g. on Schur
$\sigma$-groups).

In \cite{S5}, Scholz returns to the topic of the $(S-1)$-th power
residues of units.  Let $S$ be a generator of the Galois group of the
cyclic extension $K/\Q$ of prime degree $\ell$, and set $X = S-1$.
We recall that an element $\alpha \in K^\times$ coprime to $\ell$ is
an $(\ell, X^\lambda)$-th power residue (for $1 \le \lambda \le \ell$)
modulo a prime ideal $\fp$ if
$$ \alpha \equiv \xi^{X^\lambda} \eta^\ell \pmod {\fp} $$
for some $\xi, \eta \in K^\times$; we agree to write
$[\alpha,\fp]_\lambda = 1$ in this case, and $[\alpha,\fp]_\lambda \ne 1$
otherwise.

For $\lambda = \ell$, the element $X^\lambda = (S-1)^\ell = \ell
F_1(S)$ in the group ring $\Z[G]$ is a multiple of $\ell$. Therefore,
$X^\ell$-th power residues are always also $\ell$-th power residues,
and conversely. Hence, Scholz's notion of $X^\lambda$-th power
residues provides a refinement of the division into $\ell$-th power
residues and non-residues.

In the case $\ell = 3$, the following possibilities arise:
\begin{enumerate}
\item $[\alpha,p]_1 \ne 1$;
\item $[\alpha,p]_1 = 1$, $[\alpha,p]_2 \ne 1$;
\item $[\alpha,p]_2 = 1$, $[\alpha,p]_3 \ne 1$;
\item $[\alpha/p]_3 = 1$.
\end{enumerate}

Units with norm $+1$ are $X$-th powers by Hilbert's Theorem 90;
consequently, $[\eps,\fp]_1 = 1$ for all units $\eps$ with norm $1$.
Furthermore, Scholz's local-global principle holds:

\begin{prop}
  If $[\eps,\fp]_2 = 1$ for a unit $\eps$ and for all prime ideals $\fp$,
  then $\pm \eps = \eta^{X-1}$ for some unit $\eta$.
\end{prop}

\noindent

Indeed: if $[\eps,\fp]_2 = 1$, then $\eps^{S-1}$ is a cube modulo
every prime ideal, hence a cube of a unit: $\eps^{S-1} = \eta_1^3$.
Since $3 \sim (S-1)^2$, it follows that $\eps^{S-1} = \eta_1^{(S-1)^2}$,
and thus $\eps = \eta^{S-1} \beta$ for some $\beta$ with $\beta^{S-1} = 1$.
This implies $\beta \in \Q$, and since $\eps$
and $\eta$ are units, we must have $\beta = \pm 1$.

The question of whether $[\eps,\fp]_2 = 1$ or not can be decided by
computing the cubic residue symbols modulo the three prime ideals
$\fp$ lying above $p$:
\begin{itemize}
\item If $[\eps,p]_3 = 1$, then $(\eps,\fp)_3 = 1$;
\item If $[\eps,p]_3 \ne 1$ but $[\eps,p]_2 = 1$, then $(\eps/\fp)_3 =
  (\eps/\fp')_3 = (\eps/\fp'')_3$ have the same cubic residue
  character with respect to all three prime ideals above $p$, and
  these are pairwise distinct precisely when $[\eps,p]_2 \ne 1$.
\end{itemize}
 
\medskip\noindent
Furthermore, Scholz proves the following reciprocity law:

\begin{thm}
  Let $\ell$ and $p \equiv q \equiv 1 \pmod{\ell}$ be odd primes, and
  let $\eps_p$ and $\eps_q$ be the fundamental units of the subfields
  $K_p^\ell \subset \Q(\zeta_p)$ and $K_q^\ell \subset \Q(\zeta_q)$ of
  degree $\ell$, respectively. Then $[\eps_p,q]_2 = 1$ if and only if
  $[\eps_q,p]_2 = 1$ as well.

  In the case $\ell = 2$, the relation $[\eps_p,q]_2 = [\eps_q,p]_2$
  also holds provided that $p \equiv q \equiv 1 \pmod{4}$; in this
  case the quadratic subfields of the cyclotomic fields are real and
  possess a fundamental unit.
\end{thm}

\noindent
The part of the reciprocity law that concerns the case $\ell = 2$ had
already been discovered by Schönemann \cite{Schoene}.

The power-residue characters $[\eps,p]_\ell$ are closely related to
the divisibility of the class number of $K_p^\ell K_q^\ell$ by powers
of $\ell$. On these questions, the last word has not yet been spoken
(see e.g. Naito \cite{Naito}).

\section*{1930}

Olga Taussky, a student of Furtw\"angler who had studied the capitulation
of ideal classes in cyclic unramified extensions, had written a letter to
Scholz in November 1929 asking him for advice concerning group theoretic
problems.

At the meeting of the DMV in K\"onigsberg (Prussia, now Kaliningrad, Russia)
in September 1930, Arnold Scholz met Olga Taussky, whom he had met before
when he was in Vienna. Scholz and Taussky decide to investigate the
capitulation of ideal classes in cubic unramified extensions of quadratic
number fields.

In \cite{S6}, Scholz corrects an oversight in Furtwängler's work
\cite{Fuapp} on the approximation of real numbers by
rationals. Furtwängler started from a classical result of Dirichlet,
according to which, for every real number $\alpha$, there exist
infinitely many integers $x$, $y$ such that
$$ \Big| \frac{x}{y} - \alpha \Big| < \frac{1}{y^2}. $$
Hurwitz showed that the bound on the right-hand side can be improved
to $\frac{1}{\sqrt{5}\, y^2}$, and that this bound is optimal. Already
here one can recognize the influence of minimal discriminants: the
constant $\frac{1}{\sqrt{5}}$ in Hurwitz's theorem is related to the
fact that $D=5$ is the minimal discriminant of a real quadratic number
field. If one restricts to $\alpha \in \R \setminus \Q(\sqrt{5}\,)$,
this bound can even be improved to $\frac{1}{\sqrt{8}}$, and here too
$D=8$ is the smallest discriminant of a real quadratic number field
except for $\Q(\sqrt{5}\,)$.  

Kronecker \cite{KronDA} generalized Dirichlet's theorem to higher
dimensions: given real numbers $\alpha_1, \ldots, \alpha_n$, there
exist rational numbers with denominator $q$ such that
$$ \Big| \alpha_\nu - \frac{p_\nu}{q} \Big| < \frac{\delta}{q \sqrt[n]{q}} $$
holds for $\nu = 1, 2, \ldots, n$ and $\delta = 1$.
Minkowski \cite[p.~112]{MinDA} succeeded in reducing this bound to
$\delta = \frac{n}{n+1}$.

In the other direction, Borel \cite{BorDA} was able to show, by means
of a set-theoretic existence theorem, that Kronecker's approximation
theorem cannot be improved to
$$ \Big| \alpha_\nu - \frac{p_\nu}{q} \Big| < \frac{C}{q^{1+s}} $$
with $s > \frac{1}{n}$ and any constant $C > 0$.

Perron \cite{PerronDA} in turn was able to generalize Borel's result:
let $\alpha_1, \ldots, \alpha_n$ be algebraic integers in a real
number field $K$ of degree $n+1$; we assume that there is no relation
of the form
$$ \sum_{\nu=1}^n k_\nu \alpha_\nu - l = 0 $$
with rational $k_1, \ldots, k_n, l$ except for the trivial equation
with $k_\nu = l = 0$. Furthermore, let $\alpha_\nu^{(\mu)}$
($\mu = 1, \ldots, n+1$) denote the conjugates of $\alpha_\nu$; we set
$$ \rho_\mu = \sum_{\nu=1}^n |\alpha_\nu^{(\mu)} - \alpha_\nu|
    \qquad (\mu = 1, 2, \ldots, n), $$
as well as
$$ \sigma = \prod_{\mu=1}^n \rho_\mu. $$
Then the system of $n$ inequalities
$$ \bigg| \alpha_\nu - \frac{p_\nu}{q} \bigg| <
    \frac{\delta}{q \sqrt[n]{q}} \qquad (\nu = 1, 2, \ldots, n), $$
where $p_\nu$ and $q > 0$ denote rational numbers, cannot hold for
infinitely many values of $q$ as soon as $\delta < \frac{1}{n\sigma}$.

Perron handed over the topic of the analogous question for complex
quadratic number fields to his student Hasan Aral, who wrote his doctoral
thesis \cite{Aral} on it; incidentally, the example of the unramified
cubic extension of the quadratic number field $\Q(\sqrt{-182}\,)$
defined by the equation\footnote{Scholz probably constructed the polynomial
  $f(x) = x^3 + ax +b$ with the solution $a = 17$ and $b = 2\sqrt{-182}$ 
  of the equation $\disc f = -4a^3 - 27b^2 = 4$. The standard method
  of constructing unramified cubic extensions would be using ideas
  from Scholz's reflection theorem, namely solving the dipohantine equation
  $x^2 - 27 \cdot 182y^2 = a^3$. Such a solution comes from the
  fundamental unit $\eps = 701 + 30\sqrt{3 \cdot 182}$ and provides us
  with the polynomial $g(x) = x^3 - 3x - 1402$, whose roots generate a cyclic
cubic unramified extension of $k = \Q(\sqrt{-182}\,)$.}
$$ x^3 + 17x + 2\sqrt{-182} = 0 $$
in this dissertation was contributed by Scholz (Neiß \cite[p.~47]{NeissH}
also received from Scholz the example of the quadratic number field
$\Q(\sqrt{-19\,677}\,)$ with $3$-class group of rank $2$).

Hasan Aral was born on 25 March 1913 in Konya, Turkey; after his
Abitur, he received a scholarship from the Turkish government which
enabled him to learn the German language from 1932 to 1933 at the
Gymnasium Lüben (formerly Silesia, today Lubin; a website
\emph{Türkische Gastschüler am Lübener Gymnasium} has been set up
about this gymnasium, from which this information is taken) near
Liegnitz in Silesia.

In the winter semester 1934/35 he enrolled at Göttingen to study
mathematics, and after two semesters he transferred to Munich, where
he received his doctorate under Perron on 17 July 1939. Aral's
dissertation was cited by Hofreiter and Koksma, but was later improved
using Minkowski's methods.

In \cite{Fuapp}, Furtwängler was able to improve Perron's bound:

\begin{thm}
  Let $k < |D|^{-1/2(n-1)}$ be a positive real number, where $D$
  denotes the minimal discriminant of a real number field of degree
  $n$. Then there exist $n-1$ real irrational numbers $\alpha_1,
  \ldots, \alpha_{n-1}$, linearly independent over the rationals, for
  which the $n-1$ inequalities
  $$ \Big| \frac{x_i}{x_n} - \alpha_i \Big| < k |x_n|^{-1-\frac{1}{n-1}}
    \qquad (i = 1, 2, \ldots, n-1) $$
  have at most finitely many solutions in coprime integers
  $x_1, x_2, \ldots, x_n$.
\end{thm}

\noindent
Furtwängler notes that the minimal discriminant of cubic fields is
$-23$ (which he proved in his dissertation \cite{FurtD}), and that his
student J.~Mayer \cite{Mayer} studied the case of number fields of
degree $4$ and proved that the minimal discriminant here is $-275$.

To obtain asymptotic statements, Perron and Furtwängler used the family
of number fields $\Q(\sqrt[n]{2}\,)$ with discriminant $\le 2^{n-1} n^n$.
In fact, the discriminant of these fields is almost always exactly
$D = 2^{n-1}n^n$; for prime degrees $n = p$, the discriminant is
genuinely smaller precisely when $2^{p-1} \equiv 1 \pmod{p^2}$, which
first occurs for $p = 1093$, where $\Delta = 2^{p-1}p^{p-2}$.

Scholz later returned to the topic of minimal discriminants once more
in \cite{S21}. Not all of the conjectures stated in Furtwängler's
paper have been confirmed in the form formulated there; see Langmayr
\cite{Langm1,Langm2}.

Further generalizations of these questions can be found in numerous
papers by Cassels and Davenport on simultaneous Diophantine
approximation from the 1950s.

The paper \cite{S7} deals with the structure of ideal class groups and
unit groups of abelian fields, especially those that are composita of
two cyclic extensions of prime degree $\ell$.

To this end, let $K_p^\ell$ denote the subfield of degree $\ell$ of
the field of $p$-th roots of unity, where $p \equiv 1 \pmod{\ell}$ is
prime (or $p = \ell^2$). As the simplest example, Scholz investigates
the composita $L = K_p^\ell K_q^\ell$.

For $\ell = 2$, the situation is very clear: if $\eps_m$ denotes the
fundamental unit of the quadratic number field $\Q(\sqrt{m}\,)$, then
the unit group $E_L$ of $L$ is given by
$$ E_L = \begin{cases} 
         \la -1, \eps_p, \eps_q, \sqrt{\eps_{pq}} \ra 
            & \text{ if }  N\eps_{pq} = +1, \\
         \la -1, \eps_p, \eps_q, \sqrt{\eps_p\eps_q\eps_{pq}} \ra 
            & \text{ if }  N\eps_{pq} = -1. \end{cases} $$

Scholz calls a subfield $K$ of $L$ a \emph{norm field} if the norm map
from $L$ to $K$ is surjective on the unit group: $E_K = N E_L$.

Thus, either exactly one subfield is a norm field (namely
$\Q(\sqrt{pq}\,)$ if $N\eps_{pq} = +1$), or all three are norm fields
(namely when $N\eps_{pq} = -1$). Scholz will examine this case in
greater detail in \cite{S15}.

If $\ell$ is an odd prime, let $F$ be the Hilbert
$\ell$-class field of $L$, that is, the maximal unramified abelian
$\ell$-extension of $L$. The Galois group $G = \Gal(F/L)$ is
isomorphic to the $\ell$-class group of $L$ and, by Chebotarev, is
generated by generators of the inertia groups, hence by elements $S_1$
and $S_2$ with $S_1^\ell = S_2^\ell = 1$. The commutator subgroup $G'$
is symbolically generated by the commutator $A = S_2^{-1} S_1^{-1} S_2
S_1$, that is, $G' = \langle A^{F(S_1,S_2)} \rangle$.  

Let $\fM = \Ann(A)$ be the annihilator of $A$ in the group ring
$\Z[G]$, that is, the additive group of all elements $F$ in the group
ring such that $A^F = 1$. One easily verifies that some power of $X =
S_1 - 1$ and some power of $Y = S_2 - 1$ lie in $\fM$.

The subfields $K_p^\ell$ and $K_q^\ell$ of $L$ are called stem
fields, while the other $\ell-1$ subfields $K_\lambda$ of degree $\ell$
are termed intermediate fields by Scholz, a term which he later
replaces in \cite{S14} by middle fields.

If the norm map from $L$ to $K_\lambda$ is surjective on the unit
groups, i.e., every unit in $K_\lambda$ is the norm of a unit from
$L$, then he calls $K_\lambda$ a norm field. Scholz shows that
there are only three possibilities: 
\begin{enumerate}
\item All $\ell+1$ subfields are norm fields;
\item Exactly $\ell$ subfields are norm fields;
\item No subfield is a norm field.
\end{enumerate}

If the $\ell$-class group of $K$ is elementary abelian, Scholz can
establish further connections with $S-1$-th power residues of units.

\section*{1931}

In 1931, Lothar Heffter retires. Heffter was born in 1862 in Köslin
(now Koszalin, Poland), studied in Heidelberg and Berlin, and received
his doctorate there under Lazarus Fuchs. This was followed by
positions in Gießen (where he habilitated), Bonn, Aachen, and Kiel,
before he was appointed full professor in Freiburg in 1911. The
assistant position held by Scholz is now assigned to Doetsch.  Scholz
publishes a short paper \cite{S8} on ``Zermelo's new theory of sets'' and
simplifies his earlier considerations on disposition groups in \cite{S9}.
In \cite{S10} Scholz studies ``limitation theorems''
which essentially state that the ideal groups associated
to an extension $ K/k $ in the sense of class field theory are already
determined by the maximal abelian subextension contained in $ K/k $.
More precisely: For an extension $ F/k $, let $ \Tak_\fm(F/k) $ denote
the group of all ideals in $ k $ coprime to $ \fm $ that are norms of
ideals from $ F $. If $ L/k $ is an extension whose maximal abelian
subextension is $ K/k $, and if $ \fm $ is a multiple of the conductor
of $ K/k $, then $ \Tak_\fm(L/k) = \Tak_\fm(K/k) $.  For a modern
presentation of this fact, see \cite{GrLe}.

In March, Scholz sends Hasse a counterexample to the conjecture
expressed by Hasse, demonstrating that the norm theorem (if $K/k$ is
abelian, then an element of $k$ that is a local norm at every place is
already a global norm) is false. He recognizes that the validity of
this norm theorem would imply that $p$-class field towers terminate
after the first step; Hasse thereupon constructs his counterexample
$\mathbb{Q}(\sqrt{13},\sqrt{-3})$, which he publishes in~\cite{HasseNm}.
In 1934, Witt also found a counterexample in the case of function
fields.

\section*{1932}

In the spring of 1932, Scholz, together with Zermelo, participates in
the ``Hellas Tour'' for teachers and students of German Gymnasien;
aboard the steamship {\em Oceana}, they travel mainly in Greece from
April 19 to May 4, but also visit places in Italy and North Africa.

In 1932, \cite{S11} appears---a smaller paper, but today one of
Scholz's best-known publications. There, Scholz shows using the
simplest class field-theoretic means that the divisibility by $3$ of
the class numbers of the quadratic number fields
$\mathbb{Q}(\sqrt{m}\,)$ and $\mathbb{Q}(\sqrt{-3m}\,)$ are coupled in
such a way that the $3$-ranks of their class groups differ by at most $1$.
The deeper reason for this at first glance surprising phenomenon is that
the cubic Kummer extensions of $K = \mathbb{Q}(\sqrt{m},\sqrt{-3}\,)$,
which are abelian and unramified over $F = \mathbb{Q}(\sqrt{m}\,)$, are
of the form $K(\sqrt[3]{\alpha})$ for an
$\alpha \in k = \mathbb{Q}(\sqrt{-3m}\,)$, where $(\alpha)$ is a cube
of an ideal in $k$. Thus, unramified cyclic cubic extensions of $F$ are
related to ideal cubes in $k$, and the deviation in the $3$-rank of the
two class groups is connected with the existence of units, since under certain
conditions $\alpha$ can also be a unit.

This ``Scholz reflection theorem'' was later generalized by
Leopoldt~\cite{LeoSp} and has since been the subject of numerous further
investigations, for example by Oriat~\cite{Ori76,Ori79} and
Satg\'e~\cite{Sat76,OS79}, as well as Brinkhuis~\cite{Bri95}.

In~\cite{S12}, Scholz investigates ideal classes and units in non-cyclic
cubic number fields $K_3$ and in their normal closure $K_6$, that is, the
compositum $K_6 = K_2 K_3$ of $K_3$ and the quadratic subfield
$K_2 = \mathbb{Q}(\sqrt{d}\,)$, where $d = \disc K_3$ is the discriminant of
the cubic field. Starting from the analytic class number formula, Scholz
shows that the relation
$h_6 R_6 = h_2 h_3^2 R_2 R_3^2$
holds for the class numbers $h_j$ and regulators $R_j$ of the subfields $K_j$.
By a careful examination of the behavior of the units, Scholz finds that
$h_6 = \frac{1}{3^m} h_2 h_3^2$
holds for $m \in \{1, 2, 3\}$. Scholz obtains strong restrictions on the
structure of the Galois group of the $3$-class field of $K_6$. Parts of
Scholz's results were rediscovered by F.~Gerth~\cite{Ger3,Gerth}.

\section*{1933}

In mid-March, Scholz submits his paper~\cite{S14} to the
{\em Mathematische Annalen}, presumably via Emmy Noether. In it, Scholz
continues his investigations into the construction of number field
extensions with prescribed Galois groups, now turning, after having
already treated disposition groups, to the much more difficult problem
of constructing extensions whose Galois groups are {\em branch groups}.

He finds a necessary condition (the principal irreality condition,
Hauptirrealitätsbedingung) for a compositum of cyclic extensions with
pairwise coprime discriminants to be extendable to a branch field
(Zweigkörper). At the end of the paper, he announces that he can now
prove the existence of fields with arbitrary two-step groups generated
by two elements. An increasingly close connection between embedding
problems and norm residues emerges, which Scholz develops over the
following years into his theory of knots.

In the Reichstag election of 5 March, the National Socialists gain a
majority.  On 24 March, the Enabling Act is passed, effectively
stripping parliament of its power.  This is followed on 7 April by the
``Law for the Restoration of the Professional Civil Service'', which
forces the immediate retirement of all non-Aryan civil servants.

Gustav Doetsch\footnote{See Remmert \cite{RemmDS}} (1892--1977) had
served in the First World War; afterwards he joined the peace movement
(Friedensbund Deutscher Katholiken and Deutsche Friedensgesellschaft)
and published letters against rearmament and militarism; in addition
he had signed a petition for Emil Julius Gumbel in 1931, who had run into
trouble with nationalistic students because of his pacifism. As soon
as the Nazi's came to power, Doetsch switched sides.

In April, Loewy is dismissed from his position because he is Jewish;
Doetsch thereby becomes the sole full professor in Freiburg.  He
arranges the succession in line with his own views and states that
algebra and number theory -- like all fields under Jewish influence
-- have become entirely abstract.

Just three months after the National Socialists' seizure of power,
Scholz tries to leave Freiburg. In Berlin, his doctoral advisor Issai
Schur had already been dismissed. Erhard Schmidt managed to arrange
for Schur to continue participating in a few lecture courses during
the winter semester 1933/34.  During a visit to Berlin at
Whitsun/Pentecost 1933, Scholz discusses with Schmidt the possibility
of transferring his habilitation to Berlin.

In October 1933, Hermann Weyl resigns from his professorship in
Göttingen, stating that due to the Jewish descent of his wife he feels
``out of place'' (\cite[p.~43]{MT}), and writes in his letter of
resignation:
\begin{quote}
  {\em I cannot do otherwise than wish that the new paths which the
    present government has taken may lead the German people to
    recovery and ascent. Due to the (in my conviction) unfortunate
    entanglement with antisemitism, it is personally denied to me to
    lend a hand directly and in Germany itself.}
\end{quote}

Weyl's successor in Göttingen is Hasse.  Scholz applies for a
scholarship from the Notgemeinschaft (Emergency Association of German
Science).

\section*{1934}
Scholz gives up his assistant position with Doetsch; his successor 
is Eugen Schlotter. Scholz begins to look for opportunities 
outside of Freiburg. On 24 March 1934 he writes to Olga Taussky:
\begin{quote}
  {\em Should I ever get into difficulties, Miss Noether would also
  invite me to America at some point, which incidentally suits me
  very little; but she had already offered it to me in 1932.
  Sweden, on the other hand, would suit me very much.}
\end{quote}

In connection with his application for a research fellowship from 
the Notgemeinschaft, Doetsch, as Scholz's superior, is asked 
for his opinion.

Olga Taussky herself decides to go to Emmy Noether in the USA at Bryn
Mawr for one year; on 4 June 1934 Scholz writes to her:
\begin{quote}
  {\em How wonderful that your trip to America is now definite, and that
    you are expected at Bryn Mawr. I sincerely rejoice with you and
    am all the more mentally involved, since I now also intend to
    ask Emmy Noether about a fellowship. Perhaps I'll do it already today.}
\end{quote}

In the same letter, Scholz explains the reason for his pessimistic 
prospects in Freiburg:
\begin{quote}
  {\em The matter has been set in motion because the successor to
    Loewy has suddenly appeared here, Herr Süss from Greifswald, a
    very nice person, whom I also know through Feigl, but: No
    algebraist. Now it seems that the Doetsch lecture plan, which
    practically eliminates higher algebraic lectures, has taken hold,
    so that I now feel superfluous here.}
\end{quote}

His first idea was to ask Hasse whether he could come to him in Göttingen:
\begin{quote}
  {\em I therefore immediately wrote to Hasse about transferring my
    habilitation to Göttingen, which would now be the most suitable
    thing for me and which I had already been considering for some
    time. Flat refusal: he considers `my person' not suitable for this
    `highly charged ground' now. He may well be right; perhaps he
    wants to protect me from the company of Mister W. and the like.}
\end{quote}
In this case, Scholz was probably quite naive; even with Hasse's
support, he would not have obtained a position there, and the
right-wing student body in Göttingen would have torn him to pieces in
the air.

``Mister W.'' is most likely Werner Weber, who regarded even the
German-national Hasse only as a transitional solution until the
leadership of the Göttingen Mathematical Institute would be handed over to
a proper National Socialist. Weber had even refused to hand over the keys
of the institute to Hasse on 29 May 1934. Hasse's struggle against the
hardliners Weber and Tornier took forms in the following years that
led him on 2 April 1936 to ask Vahlen for his leave of absence and
transfer. Vahlen, however, instead recalled Tornier from Göttingen to
Berlin. After Tornier's and Weber's departure, Teichmüller also no
longer felt comfortable there and likewise left Göttingen for Berlin.

One year later, when Scholz in Kiel tries to arrange a lecture slot
for the recently dismissed Zermelo there, he has to cancel his lecture
a few days beforehand because he fears that Zermelo would be badly
treated in Kiel.

The idea of securing a position for Scholz in Kiel probably goes back
to F.K.~Schmidt: in his letter to Hasse, Scholz writes on 06.07.1934

\begin{quote}
  {\em May I ask you still to [\ldots] convey to F.\,K.\ Schmidt, if he
    is still there, my heartfelt thanks that he thought of me in Kiel!} 
\end{quote}

Friedrich Karl Schmidt was born on 22 September 1901 as the son of the
merchant Carl Schmidt and his wife Elisabeth, née Vehling, in
Düsseldorf. From 1911 he attended the Reform-Realgymnasium, where he
passed his Abitur in 1920, after which he began studying theology in
Freiburg. Because as a Realgymnasium pupil he had no knowledge of
Greek, he devoted himself exclusively to Greek in the first semester
and also attended lectures on analytic geometry with {Heffter}; by the
end of the first semester he had been converted to mathematics and
studied physics and philosophy as minors. In the autumn of 1922 his
father fell ill: F.K. Schmidt therefore dropped out of university and
began a commercial apprenticeship. When his father's health improved,
F.K. Schmidt resumed his studies and attended lectures by {Loewy} and
{Krull}; however, his father died as early as 1924. Schmidt received
his doctorate in 1925 under Loewy with ``Allgemeine Körper im Gebiet
der höheren Kongruenzen'', a generalisation of the arithmetic part of
Artin's dissertation. In the following years he played a major role in
the development of class field theory for function fields.

In the autumn of 1933 Hermann {Weyl} brought him to Göttingen so that
he could take over the algebra lectures previously given by Emmy
Noether. However, his close contacts with ``Jewish colleagues'', in
particular with Courant, soon forced him to leave. In October 1934
he received a call to Jena; in 1941 Schmidt was dismissed on the basis
of an expert opinion from the state examination committee, and he
decided to give up his teaching activities and worked during the war
at the Deutsche Versuchsanstalt für Segelflug in Ainring near Bad
Reichenhall. He returned briefly to Jena in 1945 and was appointed to
Münster in 1946. In 1952 he moved to Heidelberg, where he taught until
his retirement in 1966. Schmidt died on 25 January 1977 in Heidelberg.
\medskip

The Scholz alternative plan to obtain a position in Kiel is torpedoed
by Doetsch as soon as the matter comes to his ears. On 31 July 1934 he
writes to Tornier in Göttingen, who is to pass on his judgment to
Kaluza in Kiel:

\begin{quote}
  {\em I can only say one thing, that such a lack of sense of duty has
    not yet come to my attention at German universities. [\ldots] I do
    not demand of anyone that he be a Nat[ional]-Soc[ialist], I am not
    even a party member myself. In line with the intentions of the
    government, however, I would only promote someone who at least has
    a positive attitude toward the present state. There can be no
    question of that in the case of Sch[olz]. On the contrary, here
    together with his only friend, whom he has, namely Zermelo, he is
    known as the typical grumbler who would like nothing better than
    for the nat[ional]-soc[ialist] regime---which, in contrast to the
    previous one, threatens to get tough with characters like him---to
    disappear as soon as possible.

    It is of course out of the question that such a wimp\footnote{Doetsch
      used the word ``Waschlappen''.} as Sch[olz] could exert any
    political influence.  [\ldots] But if people
    like Sch[olz] see that despite all the government's announcements
    they are still doing very well today\footnote{[Footnote by Doetsch]
      Half a year ago he secured a fellowship from the Notgemeinsch[aft],
      the matter  was initiated behind my back.} one should not be surprised
    if they continue their mischief quietly and all the more    
    brazenly. [\ldots] You will now perhaps understand that I was
    utterly astonished that this man, precisely at a time when the
    fiercest battle is being declared against such drone existences,
    is to be honored with being called from outside for a teaching
    assignment.  [\ldots] The students are completely opposed to him,
    they even wanted to take action against him in the previous
    semester, admittedly mainly because he exclusively ran around with
    Jews and Communists.}
\end{quote}

Wilhelm Magnus found this letter in the archives after the war and
immediately sent it -- of all people! -- to Wilhelm Süss\footnote{Süss
  did not get along very well with Doetsch, although both were on the
  same side of the political divide. After the war, Süss used Doetsch
  as a scapegoat, with the intention of whitewashing himself.} with the
following words:

\begin{quote}
  {\em [I enclose] a letter, more precisely a copy of a letter, which
    Doetsch wrote to Tornier in 1934. [\ldots] I would not have had to
    be so well acquainted with Arnold Scholz to feel the urgent desire
    that it be made more difficult in the future for the writers of
    such letters to slander others in this manner. Scholz always
    assured me that he could not have endured life in Freiburg much
    longer; I now see that he was right to a degree that I could not
    properly imagine at the time!}
\end{quote}

\begin{figure}[p!]
\noindent
\includegraphics[width=11.5cm]{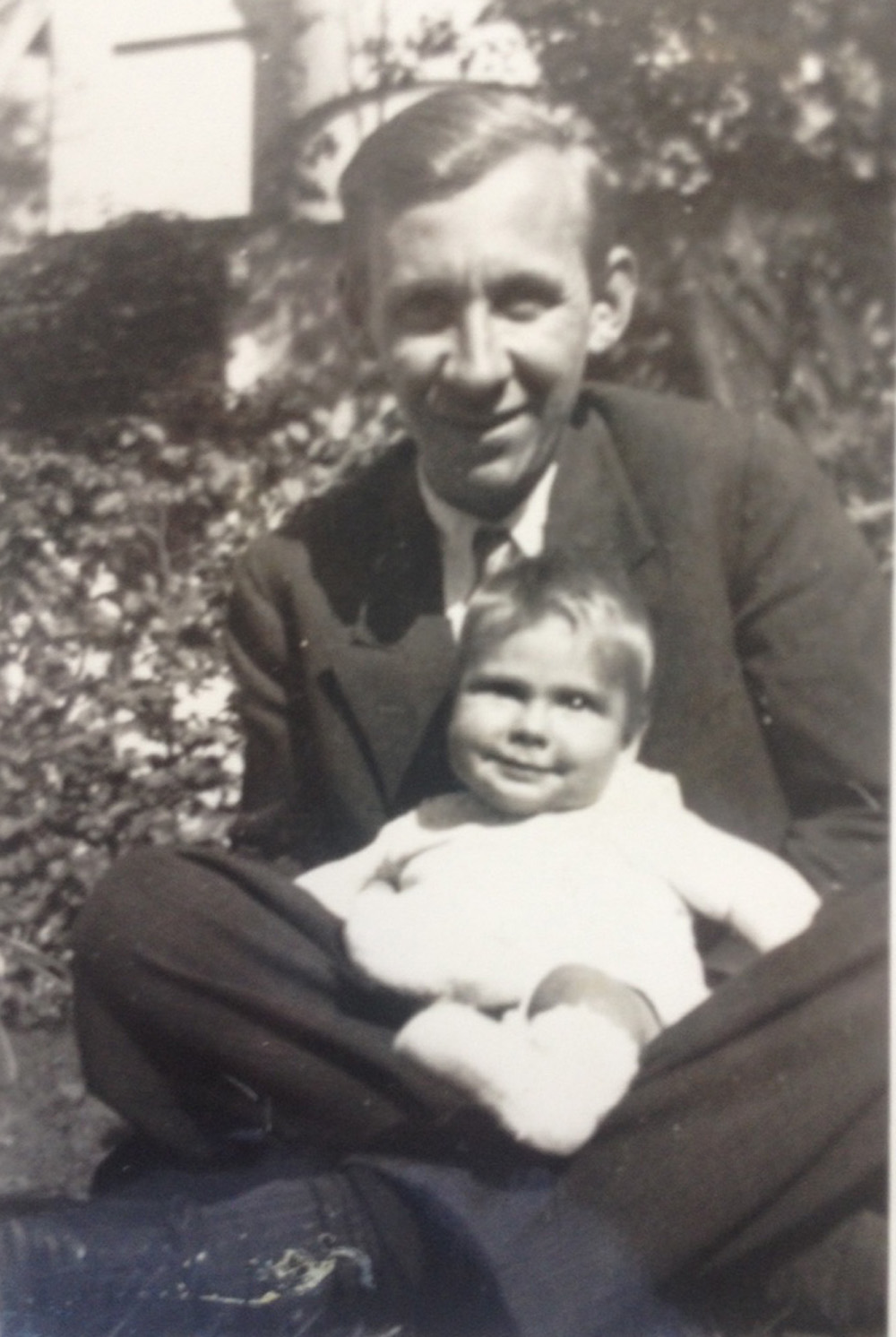}
\caption{Arnold Scholz with his nephew Nikola Korrodi}\label{Sch2}
\end{figure}

A relatively frequently cited work by Scholz is \cite{S15}, in which
he applies class field theory to investigate under what conditions
on the square-free integer $d > 0$ the equation $t^2 - d u^2 = -4$
is solvable.

Since this is a central example for Scholz's theory of knots, let us
give a few more details in the special case where $k$ is a real quadratic
number field with discriminant $d = pq$ for primes
$p \equiv q \equiv 1 \bmod 4$. Let $\eps$ denote the fundamental unit of $k$.
Here there are the following cases:
$$ \begin{array}{l|ccc}
   \rsp \text{condition} & N\eps & \Cl_2(k) & \Cl_2^+(k) \\ \hline
   \rsp (p/q) = -1                          & -1 & [2] & [2] \\
   \rsp (p/q) = +1, \ (p/q)_4 = -(q/p)_4     & +1 & [2] & [4] \\
   \rsp (p/q) = +1, \ (p/q)_4 = (q/p)_4 = -1 & -1 & [4] & [4]
\end{array} $$
If $(p/q) = +1$ and $(p/q)_4 = (q/p)_4 = +1 $, then the $2$-class group
in the strict sense is cyclic of order divisible by $8$. Scholz proves
the reciprocity law $(\eps_p/q) = (p/q)_4(q/p)_4$ for the fundamental
unit $\eps_p$ in $\Q(\sqrt{p}\,)$ already discovered by Sch\"onemann;
later Scholz will finde that this quadratic residue symbol is connected
norm residues and number knots.

The second major work, which was published in 1934 after long efforts,
is the paper \cite{S17}, on which Scholz had been working for years
together with Olga Taussky. Taussky certainly belonged to
the ``Jews and Communists'' with whom Scholz, according to Doetsch and
Schlotter, surrounded himself. The paper deals with $3$-class field
towers of imaginary quadratic number fields; more precisely, Scholz
and Taussky attempt to restrict the possibilities for the Galois group
of the $3$-class field tower by computing the capitulation in such a
way that its determination becomes possible or at least the
termination of the tower can be proved. This work has left many traces
in the mathematical literature; here we content ourselves with
references to Brink \cite{Bri84}, Heider \& Schmithals \cite{HS82},
and the many recent articles by D.~Mayer. The capitulation problem is
also, like many other problems that Scholz investigated, of a
cohomological nature. In normal unramified field extensions $K/k$ with
Galois group $G$, the capitulation kernel is isomorphic to $H^1(G,E_K)$
(see the nice survey by Jaulent \cite{Jau88}); like all problems in
algebraic number theory that have to do with units, this one is also very
difficult (or, as Artin once expressed to Taussky, ``hopeless''),
especially of course in non-cyclic extensions. Even such a simple
theorem as the one stating that in abelian unramified extensions $K/k$
a subgroup of order divisible by $(K:k)$ capitulates was proved as late as
1991 by Suzuki \cite{Suzuki}.

In July Scholz explores the possibility of transferring
his habilitation to Kiel; Hammerstein had only recently been appointed
there. In August he visits Kiel for the first time, and at the end of
October he moves there permanently.

\section*{1935}
With the end of the Notgemeinschaft fellowship, Scholz once again
begins an arduous struggle for financial independence. His attempts
to obtain a paid teaching assignment in Kiel prove more difficult than
expected, presumably also due to the continuing disruptive fire from
Freiburg. In his report from May 1935, Schlotter writes (it is a sign
of the times that even useful idiots like Schlotter are consulted for
the assessment of a mathematician such as Arnold Scholz):
\begin{quote}
  {\em Politically, Scholz is not reliable. He still does not seem to have
    come to terms with National Socialism; up to his departure from
    Freiburg he was constantly associating with Jews and former
    Communists. Scholz does not want to understand National Socialism.}
\end{quote}
Schlotter, on the other hand, understood National Socialism very well:
in 1935 he became Untersturmführer and in 1938 Hauptsturmführer of the
Waffen-SS. As late as 1992 he writes:
\begin{quote}
  {\em The Waffen-SS was an elite formation of the German Wehrmacht,
    to which the Occident owes the fact that Stalin's world conquest
    did not take place. We are proud of it.}
\end{quote}

\begin{figure}[ht!]
  \begin{center}
    \includegraphics[width=8cm]{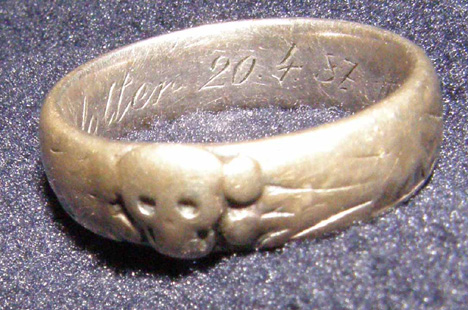}
  \end{center}  
  \caption{Schlotter's Totenkopfring}
\end{figure}

As for the elite formation: From 6 June 1944 to 17 April 1945
Schlotter commanded the Division ``Nordland'' of the
SS-Panzer-Nachrichten-Abteilung 11. Because he had recognized the
hopeless situation in Berlin and, according to his own
statements,\footnote{Sh. {\tt http://www.wehrmacht-awards.com}} did
not want to die a hero's death, he had himself transferred by a friend
to Flensburg (his successor, Hauptsturmführer Schnick, actually fell a
few days later). After the war Schlotter taught mathematics at a
secondary school until he retired in 1974---two years after the
so-called ``Radikalenerlass'' came into force in Germany, whose goal
was keeping left-leaning persons from becoming teachers, postmen or
conductors; this included in particular members of the Association of
Persecutees of the Nazi Regime and other organizations.
\medskip

From 1935 onward Scholz supervised the dissertation of G.~{Hannink} on
the displacement and non-simplicity of groups.

In the summer of 1935 it became clear that Scholz was to go to Jena as
assistant to F.K.~Schmidt. Already on 24 June F.K.~Schmidt writes to
Hasse, who had suggested Reichardt and Wecken\footnote{Franz
  Wecken passed the state examination in Göttingen in 1935 and
  received his doctorate in Marburg in 1938.} as assistants:
\begin{quote}
  {\em In addition, Blaschke has informed me in recent days of the
    difficulties with which Arnold Scholz unfortunately still has to
    deal in Kiel. According to this, Scholz has still not received
    the teaching assignment promised to him and is therefore without
    means after the expiry of his Notgemeinschaft fellowship. I would
    therefore like above all to endeavour to achieve something here
    for Scholz if possible.}
\end{quote}

Melanie Schönbeck\footnote{In {\em Kontinuität oder Wandel?}} writes
about Scholz:
\begin{quote}  
  {\em [His] lack of party membership as well as participation in
    National Socialist movements were also criticized by the Dean of
    the Kiel Philosophical Faculty, Ferdinand Weinhandl: ``He is
    neither in the SA nor in the SS, as far as I know, not even in the
    NSV.''  In contrast, according to Theodor Menzel, who had been
    Dean of the Kiel Philosophical Faculty before Weinhandl, Scholz
    was ``a man of genius, after [\ldots] Bieberbach and Tornier one
    of the best and most promising of the mathematical young
    generation.''}
\end{quote}

From Kiel, Scholz regularly attended lectures in Hamburg; in the
summer of 1935 Blaschke had paid a return visit to Kiel.  In his
letter to Hasse dated 15 July, F.K.~Schmidt once again mentions
Scholz:

\begin{quote}
  {\em I also thank you very much for mentioning the names Schröder
    and Wieland[t]. I had already heard about both of them last week
    from Schur, whom I visited in Berlin. With us, however, things now
    stand such that I have already definitively offered our assistant
    position to Scholz. Scholz, who as you know has been in Kiel since
    last October, has been held up there in the strangest
    way. Although he has already substituted for Hammerstein for
    several months, his transfer of habilitation was delayed until the
    spring of this year and the granting of the teaching assignment to
    him has not yet been pushed through to this day. Scholz is
    therefore, after the expiry of his Notgemeinschaft fellowship,
    i.e., since 1 April of this year, without any payment. I cannot
    suppress the feeling that Doetsch may have had his hand in all
    these difficulties. In any case, I know that Doetsch, in the
    autumn of last year, immediately contacted Kaluza in order to
    erase the extremely favourable impression that Scholz had
    left during his visit in Kiel. I therefore believe that a change
    of location is best for Scholz, and that he should go to a university
    where he is already known. I have not yet received an answer from
    Scholz, but assume that it will arrive in the next few days.}
\end{quote}

On 27 August 1935 Scholz writes to Taussky:
\begin{quote}
  {\em Now it is as good as decided that I will become assistant to
    F.K.~Schmidt in Jena. Not all formalities have been completed yet;
    but F.K.S. has paved the way so far.}
\end{quote}
F.K. Schmidt had moved from Göttingen in 1934, where he had been
attacked by radical students because of his dealings with Courant, to
Jena; there he also encountered problems, also because he had
published books by emigrants and Jewish authors for Springer-Verlag,
and so in 1941 he moved to the Deutsche Forschungsanstalt für
Segelflug in Ainring near Bad Reichenhall.

In October 1935 the negotiations between Scholz and Jena end;
the rector of Jena writes to Scholz on 28 October 1935 that during his
visit in Berlin he was told that it was desired for Scholz to remain
in Kiel.  F.K.~Schmidt also writes to Hasse on 30 October:
\begin{quote}
  {\em Allow me to report briefly to you on the further progress of
    the Scholz matter.  After rather protracted negotiations and after
    it had been possible to interest the Ministry of Culture in the
    matter, our rector has been informed that Berlin wishes Scholz to
    remain in Kiel, since several positions are open there.  I
    therefore assume that Scholz will receive a paid position in Kiel
    in some form for the next semester. In any case, it seems that the
    objections raised against him have finally been clarified as
    unfounded by the responsible authorities.  I am of course very
    pleased for Scholz's sake about this arrangement, even if I regret
    having to do without him now.  As we already discussed in
    Stuttgart, I have now offered the position to Herr Reichardt for
    the usual period.}
\end{quote}
On 9 December F.K.~Schmidt writes to Hasse again:
\begin{quote}
  {\em To my greatest surprise, I heard from Scholz a few days ago
    that he still has not received the teaching assignment in Kiel. I
    would like to hope that this is merely one of those delays within
    the Reich Ministry of Culture that are so frequent today. I really
    do not want to believe that the Reich Ministry of Culture's
    statement to our rector---that it is desired for Scholz to remain
    in Kiel---was merely a diplomatic maneuver by which one wanted to
    shake off further pressure in the matter from Jena's side. As long
    as Scholz still has the possibility of obtaining something in
    Kiel, I would consider it unwise if he were to leave his post and
    accept a foreign offer.}
\end{quote}
Thereupon Scholz inquires with Hasse whether Ore could obtain money
for him for a stay in Cambridge. Hasse promises to inquire with Ore,
but also writes that it would be unfavourable for Scholz's prospects
in Germany if he were to go abroad now.

In December Scholz writes that he would also go to New Haven. In his
letter to Hasse dated 19 December 1935 it becomes clear that Scholz
does not like to leave his home country:
\begin{quote}
  {\em But in the meantime I do have to concern myself with
    possibilities abroad, and I already reckon that if I go abroad, I
    will have to give up everything further here, at least for this
    era [\ldots].}
\end{quote}
Besides this possibility, Scholz also mentions others such as Amsterdam,
Cambridge and Yugoslavia. Hasse then writes to Ore on 21 December
1935:

\begin{quote}
  {\em We are at the moment somewhat worried about the shaping of the
    future of the number theorist A. Scholz, who is also known to
    you. He is currently in Kiel and has long been promised a teaching
    assignment. But since he evidently has an enemy somewhere, the
    granting of such a teaching assignment has not yet come about.}
\end{quote}
Ore's reply from January 1936 is discouraging:
\begin{quote}
  {\em I met Scholtz [sic!] in Hamburg last year and he made the best
    impression upon me. I am afraid I shall have to be rather
    discouraging about the possibilities for obtaining a position in
    America at the present time, but there might still be some
    chances. [\ldots] I should, however, not be too hopeful on the
    outcome of an application.

    The Institute of Advanced Study in Princeton announces a set of
    fellowships in the last number of the Bulletin of the American
    Mathematical Society. Scholtz may also apply for one of these
    [\ldots].}
\end{quote}

The paper \cite{S18}, in which Scholz proves that the embeddability of
the cyclotomic field $K_p^{\ell^h} K_q^{\ell^m}$ into a purely
ramified tower field is equivalent to the validity of the principal
irreality criterion, is dedicated by Scholz to his ``teacher I.~Schur
on his 60th birthday on 10 January 1935''. For this dedication, 
Scholz was denounced by Fritz Lettenmeyer (1891--1953). 
Lettenmeyer almost had been dismissed in 1934 because of
``political unreliability''. He joined the NSLB
(Nationalsozialistischer Lehrerbund) in 1933, the SA and the 
NSKK (Nationalsozialistisches Kraftfahrer-Korps) in 1934, the NSV
(Nationalsozialistische Volkswohlfahrt) in 1936, the NSdAP in 1937 and
the NSDDB (Nationalsozialistischer Deutscher Dozentenbund) in 1942 --
quite an impressive list, which apparently had no consequences for
him after the war: he remained professor in Kiel until his retirement
in 1948.

Scholz begins with the investigation of norm residues; the
corresponding examples, some of which he searched for jointly with
Artin's doctoral student Harald Nehrkorn\footnote{Harald Nehrkorn
  (1910--2006) was a student of Emil Artin and received his Ph.D.
  in 1933 with his thesis {\em On absolute ideal class groups and
    units in algebraic number fields}. After the war he worked as a
  teacher. In his lecture on class field theory from the winter semester
  1931/32, of which a transcript written after the Second World War
  has been preserved, Artin computes two examples, namely the field
  $\Q(\sqrt{-3},\sqrt{13}\,)$ of Hasse and the field
  $\Q(\sqrt{13},\sqrt{17}\,)$ ``of Nehrkorn and Scholz''.}, from
Hamburg, he later publishes in \cite{S19}.

Let $L/K$ be an abelian extension of number fields, and let $\cN$
denote the group of all $\alpha \in K^\times$ that are norm residues
in $L/K$, i.e., local norms in all completions of $L$. Scholz now
defines
$$ \begin{array}{lcl}
  \text{the number knot}  & &
    \cK_\alpha  = \cN / N_{L/K} L^\times, \\
  \text{the unit knot}  & &
    \cK_\eps  = E_K \cap \cN \Big/ E_K \cap N_{L/K} L^\times, \qquad \text{and} \\
  \text{the ideal knot}  & &
    \cK_{(\alpha)}  = \{(\alpha) : \alpha \in \cN\} \Big/
                          \{(\alpha) : \alpha \in N_{L/K} L^\times\}.
\end{array} $$

It follows from the definition that these knots are torsion groups,
and that their exponents divide the degree $(L:K)$.

Among Scholz's results is the statement that the sequence
$$ \begin{CD}
  1 @>>> \cK_\eps @>>> \cK_\alpha @>>> \cK_{(\alpha)} @>>> 1
\end{CD} $$
is exact; in his own words it reads as follows:
\begin{quote}
  {\em The unit knot is isomorphic to a subgroup of the total knot
    whose quotient group is the ideal knot. Thus for the individual group
    orders, the knot orders $k_{(\alpha)}$, $k_\eps$, $k_\alpha$, the
    relation
    $$ k_{(\alpha)} \cdot k_\eps = k_\alpha. $$
    holds.}
\end{quote}

Scholz also remarks that the knots for extensions with Galois group
$(\ell,\ell)$ are finite $\ell$-groups, and that in this case the
order of $\cK_\alpha$ divides $\ell$.  From earlier results it follows
that Rédei fields have a nontrivial number knot.  With the term Rédei
fields, Scholz denotes composita of cyclic extensions $K_1/K$ and
$K_2/K$ of prime degree, when their discriminants are coprime and when
every prime ideal that ramifies in $K_1/K$ splits in $K_2/K$ and vice
versa.

Furthermore, Scholz proves by an example that the number field
$\Q(\sqrt{p},\sqrt{q})$, where $p$ and $q$ are primes of the form
$4n+1$ with $(p/q)_4 (q/p)_4 = -1$, possesses a unit knot of order
$\# \cK_\eps = 2$.  In contrast, $\cK_\eps = 1$ when
$(p/q)_4 (q/p)_4 = +1$.

Even Andr\'e Weil pointed out in \cite{WeilA} the significance of the
discrepancy between norms and norm residues.  It was not until the
late 1960s that the theory of norm residues again came into focus for
various mathematicians.  Tate gave a cohomological approach to this
theory at the Brighton conference, thereby opening the door to modern
methods.  Garbanati, in \cite{Garb1}, was able to prove little that
had not already been established by Scholz, but in \cite{Garb2} he
drew attention to Scholz's work \cite{S25}.  Scholz's works on norm
residues were subsequently taken up again by Jehne \cite{Jehne},
Heider \cite{HeiderK}, Lorenz \cite{Lorenz}, Opolka
\cite{OpoKn,OpoKn2,OpSchur,Opolka}, and others.

Independently from this development, the arithmetic of tori was
developed, initially primarily by Takashi Ono, in which the number knot
$\nu_{K/k}$ is occasionally identified with the Tate-Shafarevich group
of the norm torus $R_{K/k}^{(1)}(\GL_1)$; see Pollio \& Rapinchuk
\cite{PoRa}.  Since no one has yet emerged to present the elementary
parts of this theory in an elementary way, this branch of number
theory has not received the dissemination it deserves.

\section*{1936}
After Scholz had been appointed a member of the scientific examination
office, he once again believes that his academic career can take an
upward turn and that he will receive the promised teaching assignment
by Easter.

In \cite{S19}, Scholz recognizes that the existence of knots is
connected with the question of the embeddability of certain field
extensions, an insight that ultimately leads to the major role played
by Schur multipliers in this problem.  In his paradigmatic example
$K = K_p^\ell K_q^\ell$, a nontrivial knot exists precisely when $K$ can
be embedded into a nonabelian extension of degree $\ell^3$ whose
nontrivial elements, in the case $\ell > 2$, all have order $\ell$,
while in the case $\ell = 2$ the Galois group is a dihedral group of
order $8$.  The Galois group of such an extension is generated by the
Frobenius substitutions $S_1$ of a prime ideal above $p$ and $S_2$ of
a prime ideal above $q$; the commutator $A = [S_1,S_2]$ corresponds to
an ideal class of order $\ell$ in $K$, the so-called {\em commutator class},
which Scholz had already introduced in \cite{S18}.  The
absolute norms of the ideals in this commutator class are principal
ideals generated by norm residues in $K/\Q$, without however being
norms themselves.

Another example from the work \cite{S19} is the following: let $k =
\Q(\sqrt{-3}\,)$, $\eta = 18 - 7 \sqrt[3]{17}$ be the fundamental unit
of $F = \Q(\sqrt[3]{17}\,)$; then $\eta = \eps^{S-1}$ in $K_1 = kF$,
where $S$ is the generating substitution of $K_1/k$, and $\eps$ and
$\eps^S$ are the two fundamental units of $k$.  The extension $K_2/k$
is cyclic of degree $3$ with prime ideal conductor $\fq$, where $\fq$
is to split completely in $K_1/k$.

For a cyclic extension $F/k$ of prime degree $\ell$, a prime ideal
$\fq$ in $k$, and an $\alpha \in F^\times$, denote by $r(\alpha,\fq)$
the maximal exponent $r \le \ell$ such that
$$ \alpha \equiv \mu^{(S-1)^r} \nu^\ell \pmod{\fq}. $$
Since units are $(S-1)$-power residues, we have $r(\eps,\fq) \ge 1$.
Scholz shows: the compositum $K = K_1 K_2$ possesses a unit knot
if and only if $r(\eps,\fq) = 1$, and an ideal knot otherwise.

Results such as this one are abundant in Scholz's works and deserve
far more attention than they have received in the past.  A very
convincing presentation of the theory of number knots was later given by
Wolfram Jehne (1926--2018) in \cite{Jehne}.  Here we content ourselves
with presenting the essential definition of the knots as given by Jehne.

Let $k$ be a number field; the embedding of $k^\times$ into the idele
group $J_k$ gives an exact sequence
$$ \begin{CD}
  1 @>>> k^\times @>>> J_k @>>> C_k @>>> 1,
\end{CD} $$
where $C_k$ denotes the idele class group of $k$.  The kernel of the
map $J_k \twoheadrightarrow I_k$ from the idele group onto the group
of fractional ideals $I_k$ is the group $U_k$ of idele units, which
gives us a second exact sequence
$$ \begin{CD}
  1 @>>> U_k @>>> J_k @>>> I_k @>>> 1.
\end{CD} $$
Finally, there is the classical sequence
$$ \begin{CD}
  1 @>>> H_k @>>> I_k @>>> \Cl_k @>>> 1,
\end{CD} $$
which defines the ideal class group $\Cl_k$ as the quotient of the
ideal group $I_k$ by the group $H_k$ of principal ideals.  These exact
sequences form the fundamental square
$$\begin{CD}
   @.       1 @.        1 @.    1       \\
        @. @VVV    @VVV          @VVV     \\
        1 @>>> E_k @>>> k^\times @>>> H_k @>>> 1 \\
            @. @VVV    @VVV          @VVV     \\
        1 @>>> U_k @>>> J_k      @>>> I_k @>>> 1 \\
            @. @VVV    @VVV          @VVV     \\
        1 @>>> \cE_k @>>> C_k      @>>> \Cl_k @>>> 1 \\
            @. @VVV    @VVV          @VVV     \\
         @.  1 @.       1 @.    1       \\
\end{CD} $$
of the algebraic number field $k$. Now Jehne considers exact sequences
\begin{equation}\label{ESJeh}
\begin{CD}
  1 @>>> A_k @>>> B_k @>>> C_k @>>> 1
\end{CD}
\end{equation}
of abelian groups like those from the fundamental square.  If $K/k$ is
normal with Galois group $G$, and if $N = N_{K/k}$ denotes the
relative norm, then an application of the snake lemma to the diagram
$$\begin{CD}
  1  @>>> A_K  @>>> B_K   @>>> C_K   @>>> 1 \\
      @. @VV{N}V    @VV{N}V    @VV{N}V     \\
  1  @>>> A_k  @>>> B_k   @>>> C_k   @>>> 1 \\
\end{CD}$$
yields the exact sequence 
$$\begin{CD} 
1  @>>> A_K[N]  @>>> B_K[N]  @>>> C_K[N]  @>{\delta}>>\\
         {} @. {} A_k/NA_K  @>>> B_k/NB_K   @>>> C_k/NC_K  @>>> 1,
\end{CD}$$
where $A[N]$ denotes the kernel of the norm map on $A$.  The
connecting homomorphism $\delta$ maps an element 
$a = \alpha A_K \in {}_N C_K$ to $\delta(a) = N_{K/k}(\alpha) N_{K/k} A_K$;
since $\alpha$ is an element of $B_K$ whose relative norm lands in $A_k$,
this is well-defined. Thus it follows that
$$ \im \delta = A_k \cap N_{K/k}B_K/N_{K/k}A_K =: [A,B]. $$
Jehne then calls $[A,B]$ the \emph{knot} associated to the sequence
(\ref{ESJeh}).  If one splits this sequence at the place $\delta$, one
obtains two exact sequences in which this knot appears:
\begin{equation}\label{knot1}
 1 \lra {}_NA_K \lra {}_NB_K \lra {}_NC_K \lra [A,B] \lra 1 
\end{equation}
\begin{equation}\label{knot2}
 1 \lra [A,B] \lra A_k/NA_K \lra B_k/NB_K \lra C_k/NC_K \lra 1.
\end{equation}
The exact sequences of the fundamental square thus yield six knots, 
of which one is trivial: $[U,J] = 1$, as can easily be shown.

Thus five knots remain, namely
\begin{enumerate}
\item the number knot
   $\nu = \nu_{K/k} = [K^\times,J_K] = k^\times \cap N J_K / N K^\times$;
\item the first unit knot
   $\omega = \omega_{K/k} = [E_K,K^\times] = E_k \cap N K^\times / N E_K$;
\item the second unit knot
   $\omega' = \omega'_{K/k} = [E_K,U_K] = E_k \cap N U_K / N E_K$;
\item the ideal knot
   $\delta = \delta_{K/k} = [H_K,I_K] = H_k \cap N I_K / N H_K$;
\item the idele knot
   $\gamma = \gamma_{K/k} = [\cE_K,C_K] = \cE_k \cap N C_K / N \cE_K$.
\end{enumerate}
With very little effort one now obtains from the main theorems of
class field theory

\begin{thm}[{\bf The Fundamental Knot Sequence}]
$$  1 \lra \omega_{K/k} \lra \omega'_{K/k} \lra
    \nu_{K/k}  \lra  \delta_{K/k} \lra \gamma_{K/k} \lra 1. $$
\end{thm}

\noindent
Scholz's unit knot $\omega^0$ is the quotient
$$ \omega'/\omega \simeq E_k \cap NU_K / E_k \cap NK^\times; $$ 
using the knot $\delta^0 := \im \,(\nu \lra \delta)$ we obtain from the
fundamental knot sequence
\begin{thm}[{\bf Scholz's Knot Sequence}]
$$ \begin{CD} 1  @>>> \omega^0_{K/k} @>>> \nu_{K/k} @>>>
                 \delta^0_{K/k} @>>> 1. \end{CD} $$
\end{thm}

\noindent
These knots allow an interpretation as Galois groups of certain 
subfields of the Hilbert class field:

\begin{thm}
We have
$$ \delta^0  \simeq \Gal(K_\cen/K_\gen), \quad 
   \delta \simeq \Gal(K_\cen/k^1K) \quad \text{and} \quad 
   \gamma  \simeq \Gal(K_\gen/k^1K). $$
\end{thm}

R\'edei fields such as $\Q(\sqrt{p},\sqrt{q}\,)$ with
$p \equiv q \equiv 1 \bmod 4$ and $(p/q) = +1$ admit an unramified
cyclic quartic extension $L/k$, hence $K_\cen$ is strictly bigger than
$K_\gen = \Q(\sqrt{p},\sqrt{q}\,)$. Thus it must have a nontrivial knot.

The middle term of Scholz's knot sequence is connected with the Schur
multiplier $\fM(G)$ of the Galois group $G = \Gal(K/k)$.  The Schur
multiplier of a finite group $G$ is defined as follows: a group
extension
$$ E: \begin{CD} 1 @>>> A @>>> \Gamma  @>>> G  @>>> 1 \end{CD} $$
is called \emph{central} if $A \subseteq Z(\Gamma)$, and a
\emph{splitting} (or split extension) if $A \subseteq \Gamma'$.  Schur
\cite{Sch07} has shown that the orders of the groups $A$ in central
splittings are bounded, and that all groups $A$ in maximal central
splittings are isomorphic to one another.  The isomorphism class of
these groups $A$ is called the Schur multiplier of $G$, and the
corresponding groups $\Gamma$ are called representation groups (or
covering groups).

It must have deeply satisfied Scholz to have uncovered the major role
played by the Schur multiplier---which his doctoral supervisor Schur
had introduced into representation theory---in class field theory.  In
his works on norm residues published in Crelle's Journal
\cite{S19,S25}, he did not cite Schur, or perhaps could not cite him;
only in the group-theoretic part \cite{S23}, which appeared in the
Monatshefte für Mathematik und Physik in Austria, do we find the
relevant citations of Schur's papers.

After the development of cohomology, it turned out that
$\fM(G) \simeq H^2(G,\Q/\Z)$.  This implies

\begin{thm}
  Let $K/k$ be a normal extension of number fields with Galois group $G$.
  Then there exists an exact sequence
  $$  \begin{CD} \coprod_\fP \what{\fM(G_\fP)}
      @>>> \what{\fM(G)} @>>> \nu_{K/k} @>>> 1 ,\end{CD}$$
  where $G_\fP$ denotes the Galois group of the local extension
  $K_\fP/k_\fp$, and where $K_\fP$ and $k_\fp$ are the completions of
  $K$ and $k$ at the prime ideal $\fP \mid \fp$.  Furthermore,
  $\widehat{A} = \Hom(A,\Q/\Z)$ denotes the Pontryagin dual group of
  $A$.
\end{thm}

\noindent
Since the Schur multiplier of cyclic groups is trivial, 
this theorem has the Hasse norm theorem as a corollary.

\bigskip\noindent
In addition to his first investigations of knots and their
connection with the construction of field extensions with prescribed
Galois groups, Scholz studies Furtwängler's work on rational
approximability and the related question of the size of minimal
discriminants of number fields with prescribed degree.  He publishes
his results in his article on minimal discriminants \cite{S21}.  In
doing so, Scholz takes up again a question from \cite{S6}: how can one
find towers of field extensions whose discriminants remain relatively
modest?  In this direction, apart from Perron's example of the fields
$\Q(\sqrt[n]{2}\,)$, nothing was known.

Among all number fields of degree $n$, let $K_n$ be a number field with
minimal discriminant (meaning the absolute value). Scholz considers
the real number $E_n = \sqrt[n]{|\disc K_n|}$,
which in the modern literature is called the root discriminant $\rd(K_n)$.
Results of Minkowski show
$$ \liminf_{n \to \infty} E_n \ge e^2, $$
and Blichfeldt improved this bound to
$ \ge e \pi$.
Perron's fields $\Q(\sqrt[n]{2}\,)$ yield $E_n < 2n$; 
from the existence of infinite class field towers, as Scholz writes, 
it would follow that
$ \liminf_{n \to \infty} E_n / n = 0. $

By means of cyclotomic fields $K/\Q$ one can only slightly improve 
Furtwängler's upper bound, because for such fields
$$ E_n \ge C \cdot \frac{n \log \log n}{\log n}. $$
for $K = \Q(\zeta_p)$ we have $\disc K = p^{p-2}$, so
$\rd(K) = p^{\frac{p-2}{p-1}}$.

By considering ray class fields over cyclotomic fields, however, 
Scholz succeeds in proving that
$$ E_n = O\Big( \frac{\log n}{\log \log n}\Big)^2. $$
Scholz's construction in fact yields unramified abelian extensions 
of suitable cyclotomic fields; see Lemmermeyer \cite{Lcyc3}.

Scholz's conjecture that $\lim_{n \to \infty} E_n / n = 0$ remains
open even despite the techniques of Golod and Shafarevich.  The
difficulty lies in proving the existence of prime-degree extensions
with small discriminant.  Similarly, Scholz's conjecture that the
minimal discriminant for fields of degree $n = 400$ is smaller than
that for $n = 397$ is still open.

There are a whole series of further questions from this area that seem
interesting.  For example, Ankeny \cite{Ank} has shown that there
exists a constant $c > 0$ such that
$\log \rd(K) > c \cdot \log_m [K:\Q]$ for all normal extensions $K/\Q$
whose Galois group $G$ is metabelian of step $m \ge 1$, for which the
$m$-th commutator subgroup $G^{(m)} = 1$ vanishes, but the $(m-1)$-th
does not.  Here, $\log_m(x) = \log \log \cdots \log(x)$ denotes the
$m$-fold iterated logarithm.

It is, to my knowledge, not known whether this bound is optimal.
For $m = 1$, i.e., abelian groups, the fields
$K = \Q(\zeta_p)$ have root discriminant $p^{\frac{p-2}{p-1}}$, and thus
$ \log \rd(K) = \frac{p-2}{p-1} \log p < \log (p-1) = \log [K:\Q]. $
Scholz's result shows that such families of fields also exist for $m = 2$, 
since for the metabelian fields he constructed the inequality
$$ \log \rd(K) < 2 \log_2 [K:\Q] $$
holds.

From de Gruyter, Scholz receives the offer to write an introduction to 
elementary number theory for the Göschen Collection.

\section*{1937}
In search of a permanent position, Scholz writes to Hasse on 3
November 1937 that he might be well placed in T\"ubingen, where a
successor to Kamke was currently being sought.  Kamke was married to a
woman of Jewish descent, although she had converted to Protestantism.
Already in April 1937, Kamke was given the choice of divorcing his
wife or losing his position; Kamke (see Segal \cite[pp.~105]{Segal})
chose the latter and was dismissed on 1 November.  That Kamke survived
the war years unharmed is due, among other things, to the fact that
S\"uss kept his protective hand over him.

In 1937, Scholz publishes his paper \cite{S20} on the construction of
algebraic number fields whose Galois group is a prescribed $p$-group.
The fact that he dedicates this article ``to the memory of Emmy
Noether'' will hardly have helped his efforts to secure a position.
In this paper, Schur's representation groups appear in connection with
embedding problems, once again demonstrating how closely such
questions are linked to the theory of knots that Scholz would later
develop.

In 1937, Scholz's Problem 253 on ``addition chains'' appears in the
Jahresberichte, which, due to its simple formulation and its
complexity, quickly became famous and notorious.  Alfred Brauer
succeeds in \cite{BraAC} in answering some of the problems, but
Scholz's actual conjecture has remained open to this day and is
discussed in Guy's ``Unsolved Problems in Number Theory''
\cite[C6]{Guy}.  It was only later discovered that the problem of
addition chains had already been formulated in 1894 by H.~Dellac in
\cite{Dellac}.

Scholz revises his second paper on norm residues for the sixth time;
this work \cite{S25} will only appear in 1940 as one of his last
publications.  Since the mathematical language suitable for describing
Scholz's insights (above all cohomology theory) has not yet been
developed, this paper is among the most difficult to read.  Here too,
Scholz has hidden many gems so well within the text that his priority
was only recognized much later: on p.~227, for example, he writes:
\begin{quote}
  {\em Moreover, it is shown that every finite group occurs as the
    Galois group of an unramified relative extension.}
\end{quote}

This immediately implies that there exist arbitrarily long
class field towers.  Artin had already proved this theorem in 1934,
but did not publish it: cf.\ \cite[15.1]{FLR} and \cite[7.27]{LRH}.
Later, Fr\"ohlich \cite{FroeUn} and Uchida \cite[Cor.~1]{Uchi}
rediscovered this result; a strengthening for $p$-groups can be found
in Nomura \cite{NomGal}, and Ozaki \cite{Ozaki} was finally able to
show that for every $p$-group $G$ there exists not only an unramified
extension $K/k$ with Galois group $G$, but moreover one may
additionally require that $K$ is the $p$-class field tower of $k$.

\section*{1938}
In October 1938, Hammerstein asks Süss -- who had been proposed as
chairman of the DMV in August 1937 and was elected shortly afterwards
at the annual meeting in Bad Kreuznach -- for an expert opinion on
Scholz (Hasse was also contacted; he replies (see
\cite[p.~62--63]{RemmFB}):
\begin{quote}
  {\em I must honestly confess that, in this sense, I would certainly
    not yet propose Mr.~Scholz to the government of the Third Reich at
    the present moment.  His attitude has left much to be desired for
    years.}
\end{quote}
In 1938, Süss also actively participated in removing Issai Schur from
the editorial board of the \textit{Mathematische Zeitschrift} (see
Remmert \cite{RemmG}). The highly ambiguous role that Süss played during
the Third Reich did not harm him in any way after 1945: Süss was certainly
a master of walking the tightrope between diplomatic skill and
opportunism. Süss faced almost no problems after the war, although
Ferdinand Springer, backed by F. K. Schmidt, accused Süss of having
tried to have Schmidt arrested in 1939 when Schmidt was visiting the
US (see Remmert \cite{RemmFB}). Springer had sent Schmidt to the USA in
connection with the Zentralblatt; Süss accused Schmidt of having strong
ties to Jewish emigrants and suggested he himself should go. After the
war, Süss accused Schmidt of having given an erroneous report.

Scholz completes his work on his \textit{Introduction to Number
  Theory} \cite{S22} in the Göschen Collection.  This book begins with the
derivation of the most important fundamental laws of the natural
numbers from the axioms\footnote{Here Zermelo's influence can be seen.},
develops elementary number theory within the
framework that has been essentially standard since Gauss, and then
provides an introduction to the theory of binary quadratic forms.  In
the final chapter, Scholz explains various tricks and methods for
computations in number theory.

Scholz's hopes for a position in Vienna are dashed; depressions make
it hard for him to do mathematics.  In the following years, Scholz
limits himself to completing works already begun.  In particular, the
continuation of the work on minimal discriminants, which had already
been announced, never appeared: On 27 July 1938, Scholz responds to an
inquiry from Hasse concerning minimal discriminants of number fields
with prescribed degree, and in doing so introduces the Dirichlet
series $D_{n,G}(s) = \sum |D_K|^{-s}$, where the summation extends
over all fields of degree $n$ whose normal closure has Galois group
$G$.  He also communicates various results on the abscissa of
convergence $a(G)$ of such series; in particular, for degree-$4$
extensions,
\begin{itemize}
\item $a(D_4) = 1$ for the dihedral group $D_4$ of order $8$, 
\item $a(V_4) = a(C_4) = \frac{1}{2}$ for the Klein four group $V_4$
  and the cyclic group $C_4$ of order $4$.
\end{itemize}
Moreover, Scholz correctly conjectures that $a(A_4) = \frac{1}{2}$ and
$a(S_4) = 1$.  These results are connected with asymptotic questions
about the number of number fields with given Galois groups, which have
been intensively studied in recent years.

One of the central conjectures in the theory of counting Galois
extensions is the following: If $N(G,x)$ denotes the number of all
extensions of degree $n$ with discriminant $\le x$ whose normal
closure has Galois group $G$, then
$$ N(G,x) = c \cdot x^{a(G)} (\log x)^{b(G)} $$ holds for certain
constants $a$, $b$, $c$ that depend only on the base field and $G$.
Writing $Z_G(s) = \sum D_K^{-s} = \sum a_n n^{-s}$, we have
$\sum_{k=1}^n a_k = N(G,n)$.  From known properties of Dirichlet
series, it then follows that $a(G)$ equals the abscissa of convergence
of $Z_G(s)$.  In particular, the values claimed by Scholz, namely
$a(D_4) = 1$ and $a(C_4) = a(V_4) = \frac{1}{2}$, are all correct.
Moreover, Scholz's conjectures $a(S_4) = 1$ and $a(A_4) = \frac{1}{2}$
turned out to be true, as has the claim that
$\lim_{s \to 1/2} \frac{Z(s)}{V(s)} = 0$ on account of $b(C_4) = 0$
and $b(V_4) = 1$. Nothing is known about Scholz's proofs.

\section*{1939}
In 1939, Scholz publishes an original introduction to number theory
\cite{S22} in the Göschen Collection.

From 26 to 30 June, Hasse organizes a group theory conference in
Göttingen, to which he has invited Philip Hall.  The first four days
each begin with a two-hour lecture by Hall, after which Magnus,
Zassenhaus, Grün, Speiser, Wielandt, Witt, van der Waerden, Krull,
and Specht give talks; Scholz delivers a lecture on 29 June from 18:00
to 19:00 on {\em Extension and splitting of groups}, thus on the
content of his papers \cite{S23,S24}.  Many of these lectures were
published in volume 182 of Crelle's Journal.

In \cite{S23}, Scholz investigates the group-theoretic side of field
constructions as they are needed in his theory of norm residues.  In
order to convey at least a rough idea of Scholz's terminology, we here
provide a brief introduction to the presentation of the theory as it
was developed by Heider \cite{HeiderK} and later by Steinke
\cite[\S~2]{Steinke}.

To this end, let $K/k$ be a normal extension of number fields. 
If an extension $L/K$ is normal over $k$, then we have an exact sequence
$$ \begin{CD}
  1 @>>> \Gal(L/K) @>>> \Gal(L/k) @>>> \Gal(K/k) @>>> 1
\end{CD} $$
A normal extension $L/K/k$ is called
\begin{itemize}
\item \emph{central over $K/k$} if $\Gal(L/K) \subseteq Z(\Gal(L/k))$,
  that is, the Galois group of $L/K$ lies in the center of the Galois
  group of $L/k$.
\item \emph{abelian over $K/k$} if $\Gal(L/K) \cap \Gal(L/k)' = 1$,
  that is, the Galois group of $L/K$ and the commutator subgroup of
  the Galois group of $L/k$ are disjoint.
\end{itemize}
Observe that if $L/K$ is abelian over $K/k$, then $L/k$ is not necessarily
abelian. For example, $L/K$ is abelian over $K/k$ whenever $L = KF$ for some
abelian extension $F/k$.

Now let $L_1/k$ and $L_2/k$ be central extensions over $K/k$.  If
there exists a central extension $L/k$ over $K/k$ that contains both
$L_1$ and $L_2$ and which is abelian over $L_1/k$ and over $L_2/k$,
then $L_1$ and $L_2$ belong to the same \emph{genus}.

Opolka has shown in \cite{OpGesch} that the genera of central
extensions of a Galois extension $K/k$ of global or local fields
correspond bijectively to the quotient groups of $H^{-1}(G,C_K)$,
where $C_K$ denotes the idele class group of $K$ in the global case
and the multiplicative group $K^\times$ in the local case.

Here once again it becomes apparent how Scholz, with his own language,
created a substitute for the cohomological structures that were not
yet available at the time.

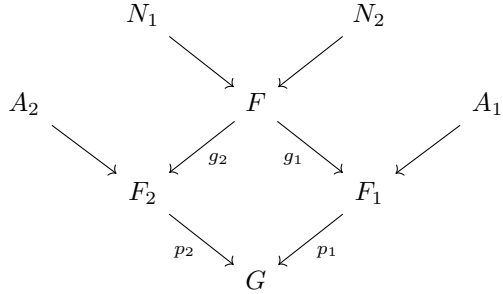
\begin{figure}[ht!]
  \begin{center}
  \begin{tikzcd}
        & N_1 \arrow[rd] &           & N_2 \arrow[ld] &            \\
  A_2 \arrow[rd] & & F \arrow[ld,"g_2"] \arrow[rd,"g_1",labels=below left]
          & & A_1 \arrow[ld] \\
    & F_2 \arrow[rd,"p_2",labels=below left] &  & F_1 \arrow[ld,"p_1"] &   \\
         &         & G           &     &               
  \end{tikzcd}
  \end{center}
\caption{Scholz's group genera}\label{SchGG}
\end{figure}

If, in addition, there exists an abelian extension $M/k$ that contains 
the maximal abelian subextension of $K/k$ and for which 
$L_1 M = L_2 M$ and $L_1 \cap M = L_2 \cap M$ hold, 
then one says that $L_1$ and $L_2$ belong to the same class, 
or that they arise from each other by \emph{abelian crossing}.
Scholz calls the group $F$ the \emph{relative product} of $F_1$ and
$F_2$ if $F$ is chosen minimally in Fig.~\ref{SchGG}.

Hilbert introduced the technique of abelian crossing in his proof of
the Kronecker--Weber theorem: there he showed that one can make an
abelian extension $K/\Q$ disappear by abelian crossing with cyclotomic
fields, from which it immediately follows that $K$ is contained in a
compositum of cyclotomic fields.  Chebotarev was able to prove his
density theorem using the same technique, and Artin employed it in the
proof of his reciprocity law.  Despite its outstanding importance,
Scholz's penetration of this concept with the help of group theory
found little resonance.

From the group-theoretic point of view, the matter comes down to the following: 
let $G = \Gal(K/k)$ be a finite group and set $F_j = \Gal(L_j/k)$; 
then the sequences
\begin{equation}\label{ExS1} 
\begin{CD} 1  @>>>  A_j @>>> F_j @>>> G @>>> 1 \end{CD} 
\end{equation}
are exact, and if we identify $A_j$ with its image in $F_j$, then
$A_j \subseteq Z(F_j)$.  If there exist abelian group extensions
$$ \begin{CD} 
   1  @>>>  N_j @>>> F @>>> F_j @>>> 1, 
   \end{CD} $$
for which the diagram in Fig.~\ref{SchGG} is commutative, 
then the exact sequences (\ref{ExS1}) belong to the same genus.

Furthermore, if there exists a normal subgroup $N$ of $F$ such that
$F/N$ is abelian and with $NN_1 = NN_2$ and $N_1 \cap N = N_2 \cap N = 1$,
then the exact sequences (\ref{ExS1}) belong to the same
class, or arise from each other by abelian crossing.

\section*{1940}
In the year 1940, G.~Hannink receives his Ph.D. under Scholz in Kiel 
with the thesis\footnote{Hannink cites Schur in a footnote.}
``Transfer and non-simplicity of groups''.  In May, Scholz is drafted,
trained as a radio operator, and sent to the East.

\section*{1941}
In the year 1941, Scholz is assigned as a mathematics teacher at the 
Naval Academy Flensburg-M\"urwik, together with Ott-Heinrich Keller
(1906--1990). On 28 February, Scholz is in Kiel to attend the funeral
of Hammerstein, 
who had died on 25 February ``due to malnutrition affecting the lungs''.

Scholz organizes a colloquium on 19 July in G\"ottingen on the occasion 
of Zermelo's 70th birthday; Zermelo himself delivers three lectures, 
and the other speakers are Konrad Knopp and Bartel van der Waerden.

In December, Scholz learns that he cannot count on remaining in Kiel 
after the end of the war.

Mathematically, he occupies himself with infinite algebraic extensions
of number fields, which had already been studied before him by Krull
and Herbrand; Moriya \cite{Mori36a,Mori36b} had also concerned himself
in 1936 with the development of a class field theory over infinite
algebraic number fields.  The corresponding paper, which Scholz had
originally intended to deliver to Hasse in Berlin at Christmas, was
only published posthumously.

\section*{1942}

Scholz dies on 1 February 1942 in Flensburg. 
In \cite[p.~141]{Bigalke} it is stated about this:
\begin{quote}
  {\em Occasionally there were also musical contacts with the number
    theorist Arnold Scholz, who taught as a lecturer at the University
    of Kiel.  Heesch later reported of him that he had hardly known
    any other person who could listen to music so well and understand
    its spirit.  Scholz considered Franz Schubert the greatest
    musician in the entire history of mankind and his String Quartet
    in C major, op.~163 as the greatest musical work.
    While Heesch was active during the war as a ``civilian teacher''
    at the Naval Gunnery School in Kiel, Scholz taught as a civilian
    teacher at the similarly named school in Flensburg.  They had
    mutual friends at whose homes they were often invited.  At one of
    these festive evenings in January 1942, where ``despite the
    darkness of the times, there was also dancing'', Scholz struck the
    hostess by his unusual thirst.  Eight days later Arnold Scholz was
    dead.  Cause: ``galloping diabetes''.}
\end{quote}
``Galloping'' was an archaic adjective used in medicine
to describe diseases that advance quickly and aggressively.

Heesch had already met Olga Taussky in G\"ottingen in 1931, who
occasionally joined him, Franz Rellich, Ernst Witt, Stefan
Cohn-Vossen, and Fritz John for communal lunches.  Scholz's mother
sends his estate to relatives in Poland; during the expulsion of the
Germans in 1945, all documents are lost.

Regarding the history of Scholz's estate, Hasse wrote to Olga
Taussky on 29 October 1949:
\begin{quote}
{\em Mrs.~Scholz\footnote{Mars. Scholz is his mother, who was living in
    Zurich at the time.}  initially wrote to me about the estate as
  follows: It was stored in a hard case and several suitcases with
  his relatives Josef Stenzel in G\"unthersdorf near
  Brieg/Oder\footnote{Today: Brzeg}.  When the Poles arrived, these
  relatives had to vacate their property, and Arnold Scholz's entire
  estate along with some of his furniture was ``taken into custody''
  by the Poles.  Given the situation, it is hardly possible for us
  today to gain access to these items.}
\end{quote}
This indeed proved to be the case, so that the entire estate of Scholz
(as well as that of Emmy Noether) must be considered completely lost;
Hasse himself calls this ``a very painful loss''.

\section*{1943}

Scholz's last paper \cite{S26}, which he dedicated to his recently
deceased colleague Hammerstein, is published posthumously.  It
concerns a new foundation for the number theory of algebraic field
extensions of infinite degree, in which he recovers earlier results of
Krull \cite{Krull1,Krull2} and Herbrand \cite{Herbinf} by different
means.

At that time, it was in no way foreseeable that these ideas would one
day bear rich fruit: following the groundbreaking work of Iwasawa,
infinite extensions of number fields became an indispensable tool
for studying arithmetic questions in finite extensions.

\section*{1949}
Ott-Heinrich Keller (1906--1990) dedicates his paper \cite{Kell}
``to the memory of Arnold Scholz''.

\section*{1952}
Olga Taussky publishes an obituary \cite{T52a} for Arnold Scholz.

\section*{1969}
Emma Lehmer discovers the reciprocity law $(\eps_p/q) = (p/q)_4(q/p)_4$
for primes $p \equiv q \equiv 1 \bmod 4$ with $(p/q) = 1$ and asks
Hasse (letter from October 24) if he knows this result. Hasse replies
that he does not. Lehmer finds the result in Scholz's article and tells
Hasse about in a letter from January 24, 1970.

\section*{1979}
Wolfram Jehne publishes his major work \cite{Jehne} 
on number-theoretic knots ``In memoriam Arnold Scholz'' 
in Crelle's Journal.

\section*{1982}
On the occasion of the 40th anniversary of his death 
and the publication of Heider and Schmithals \cite{HS82}, 
Crelle's Journal publishes a photograph of Arnold Scholz.

\section*{2016}

Lemmermeyer and Roquette publish the correspondence between Hasse,
Scholz and Taussky.

\vfill \eject

\section*{Publications of Arnold Scholz}

\subsection*{Problems in the Jahresberichte der DMV}

{\bf Problems posed:} 194, vol. {\bf 45}; 208, vol. {\bf 45};
211, vol. {\bf 46}; 212, vol. {\bf 46}; 222, vol. {\bf 46};
232, vol. {\bf 46}; 249, vol. {\bf 47}; 253, vol. {\bf 47};

\medskip \noindent
{\bf Problems solved:} 153, vol. {\bf 43};  169, vol. {\bf 44};
171, vol. {\bf 44}; 175, vol. {\bf 44};


\begin{enumerate}

\bibitem[S1]{S1} (mit H. Hasse),
  {\em Zur Klassenkörpertheorie auf Takagischer Grundlage},
  Math. Z. {\bf 29} (1928), 60--69
  %

\bibitem[S2]{S2}
  {\em Über die Bildung algebraischer Zahlkörper mit auflösbarer
    Galoisscher Grup\-pe}, Math. Z. {\bf 30} (1929), 332--356  
  %
´
\bibitem[S3]{S3}
{\em Reduktion der Konstruktion von Körpern mit zweistufiger
(metabelscher) Gruppe}, Sitz.ber. 
Heidelberger Akad. Wiss. 1929, Nr. 14, 3--15
%

\bibitem[S4]{S4}
  {\em Anwendung der Klassenkörpertheorie auf die Konstruktion von
    Körpern mit vorgeschriebener Gruppe}, Jber. DMV {\bf 38} (1929), 46
  %

\bibitem[S5]{S5}
  {\em Zwei Bemerkungen zum Klassenkörperturm},
  J. Reine Angew. Math. {\bf 161} (1929), 201--207
  %
  
\bibitem[S6]{S6}
  {\em Zur simultanen Approximation von Irrationalzahlen},
  Math. Ann. {\bf 102} (1930), 48--51
  %

\bibitem[S7]{S7}
  {\em Über das Verhältnis von Idealklassen- und Einheitengruppe
    in Abelschen Kör\-pern von Prim\-zahl\-po\-tenz\-grad},
  Sitz.ber. Heidelberger Akad. Wiss. 1930, Nr. 3, 31--55
  %

\bibitem[S8]{S8}
  {\em Zermelos neue Theorie der Mengenbereiche},
  Jber. DMV {\bf 40} (1931), 42--43
  %

\bibitem[S9]{S9}
{\em Ein Beitrag zur Theorie der Zusammensetzung endlicher Gruppen},
Math. Z. {\bf 32} (1931), 187--189 
%

\bibitem[S10]{S10}
  {\em Die Abgrenzungssätze für Kreiskörper und Klassenkörper},
  S.ber. Akad. Wiss. Berlin 1931, 417--426 
  %

\bibitem[S11]{S11}
  {\em Über die Beziehung der Klassenzahlen quadratischer Körper
    zueinander}, J. Reine Angew. Math. {\bf 166} (1932), 201--203
  %

\bibitem[S12]{S12}
{\em Idealklassen und Einheiten in kubischen Körpern},
Monatsh. Math. Phys. {\bf 40} (1933), 211--222 
%

\bibitem[S13]{S13}
{\em Die Behandlung der zweistufigen Gruppe als Operatorengruppe},
Sitzungsber. Heidelberg. Akad. Wiss. 1933, No.2, 17--22 
%

\bibitem[S14]{S14}
  {\em Die Kreisklassenkörper von Prim\-zahl\-po\-tenz\-grad und die
    Konstruktion von Körpern mit vorgegebener zweistufiger Gruppe I},
Math. Ann. {\bf 109} (1933), 161--190 
%


\bibitem[S15]{S15}
{\em Über die Lösbarkeit der Gleichung $t^2-Du^2=-4$},
Math. Z. {\bf 39} (1934), 95--111 
%

\bibitem[S16]{S16}
  {\em Berichtigung zu der Arbeit: ''Die Kreisklassenkörper von
    Prim\-zahl\-po\-tenz\-grad und die Konstruktion von Körpern mit
    vorgegebener zweistufiger Grup\-pe I''},
Math. Ann. {\bf 109} (1934), 764 
%

\bibitem[S17]{S17} (mit O. Taussky),
{\em Die Hauptideale der kubischen Klassenkörper imaginär\-qua\-dra\-ti\-scher
Zahlkörper: ihre rechnerische Bestimmung und ihr Einfluss auf den
Klassenkörperturm},
J. Reine Angew. Math. {\bf 171} (1934), 19--41 
%

\bibitem[S18]{S18}
{\em Die Kreisklassenkörper vom Primzahlpotenzgrad und die Konstruktion
von Körpern mit vorgegebener zweistufiger Gruppe II},
Math. Ann. {\bf 110} (1935), 633--649 
%


\bibitem[S19]{S19}
{\em Totale Normenreste, die keine Normen sind, als Erzeuger
nicht-abelscher Körpererweiterungen I},
J. Reine Angew. Math. {\bf 175} (1936), 100--107 
%

\bibitem[S20]{S20}
  {\em Konstruktion algebraischer Zahlkörper mit beliebiger Gruppe von
    Prim\-zahl\-po\-tenz\-ord\-nung I},
  Math. Z. {\bf 42} (1937), 161--188
  %

\bibitem[S21]{S21}
{\em Minimaldiskriminanten algebraischer Zahlkörper},
J. Reine Angew. Math. {\bf 179} (1938), 16--21 
%

\bibitem[S22]{S22}
  {\em Einführung in die Zahlentheorie},
  Sammlung G\"oschen. 1131) Berlin: de Gruyter 136 pp. (1939);
  2nd edition 136 pp. (1945);
  2nd edition revised by B. Schöneberg, 128 pp. (1955);
  5th edition revised and edited by B. Schöneberg, 128 pp. (1973)
  %

\bibitem[S23]{S23}
  {\em Abelsche Durchkreuzung},
  Monatsh. Math. Phys. {\bf 48} (1939), 340--352
  %

\bibitem[S24]{S24}
  {\em Zur Abelschen Durchkreuzung},
  J. Reine Angew. Math. {\bf 182} (1940), 216
  %

\bibitem[S25]{S25}
  {\em Totale Normenreste, die keine Normen sind, als Erzeuger
    nicht-abelscher Körpererweiterungen II},
  J. Reine Angew. Math. {\bf 182} (1940), 217--234
  %

\bibitem[S26]{S26}
  {\em Zur Idealtheorie in unendlichen algebraischen Zahlkörpern},
  J. Reine Angew. Math. {\bf 185} (1943), 113--126
  %

\bibitem[S27]{S27}
  {\em Spezielle Zahlkörper} ({\em Special number fields}),
  Bestimmt für: Enzykl. math. Wiss. 2. Aufl. I (unvollendet)
  %

\end{enumerate}

\newpage

\begin{thebibliography}{999}

\bibitem{Ank} N.C. Ankeny,
  {\em An improvement of an inequality of Minkowski},
  Proc. Nat. Acad. Sci. U.S.A. {\bf 37} (1951), 711--716
  %

\bibitem{Aral} H. Aral,
  {\em Simultane diophantische Approximationen in imaginären
    qua\-dra\-ti\-schen Zahl\-kör\-pern},  Diss. München, 1939 
  %

\bibitem{Azizi} A. Azizi, 
{\em Capitulation des $2$-classes d'id\'eaux de $\Q(\sqrt{d}, i)$},
Thesis Univ. Laval 1993 
%

\bibitem{BaBu} L. Bartholdi, M.R. Bush,
{\em Maximal unramified $3$-extensions of imaginary quadratic fields and 
     SL$_2(\Z_3)$},
J. Number Theory {\bf 124} (2007), 159--166
%

\bibitem{Bigalke} H.-G. Bigalke,
{\em Heinrich Heesch. Kristallgeometrie, Parkettierungen, 
     Vier\-far\-ben\-for\-schung},
Vita Mathematika Band 3 (E. Fellmann, Hrsg.), Birkhäuser 1988 
%

\bibitem{BLS1} E. Benjamin, F. Lemmermeyer, C. Snyder,
  {\em Imaginary quadratic Fields with cyclic $Cl_2(k^1)$},
  J. Number Theory {\bf 67} (1997), 229--245
  %

\bibitem{BLS2} E. Benjamin, F. Lemmermeyer, C. Snyder,
  {\em Real quadratic number fields with Abelian Gal($k^2/k$)},
  J. Number Theory {\bf 73} (1998), 182--194
%

\bibitem{BLS3} E. Benjamin, F. Lemmermeyer, C. Snyder,
  {\em Imaginary quadratic fields $k$ with $\Cl_2(k) \simeq (2,2^n)$
    and rank $\Cl_2(k^1) = 2$}, Pac. J. Math. {\bf 198} (2001), 15--32
  %

\bibitem{BLS4} E. Benjamin, F. Lemmermeyer, C. Snyder,
  {\em Imaginary quadratic fields $k$ with $\Cl_2(k)$ of type $(2,2,2)$},
  J. Number Theory  {\bf 103} (2003), 38--70
  %

\bibitem{BorDA} E. Borel,
  {\em Contribution \`a l'analyse arithm\'etique du continu},
  J. Math. Pures Appl. (5) {\bf 9} (1903), 329--375
  %

\bibitem{BraAC} A.T. Brauer,
  {\em On addition chains},
  Bull. Amer. Math. Soc. {\bf 45} (1939), 736--739 
  %

\bibitem{Bri84} J.~R.~Brink,
  {\em The class field tower of imaginary quadratic number fields
    of type $(3,3)$},  Ph. D. Diss. Ohio State Univ (1984),  121pp 
  %

\bibitem{Bri95} J.~Brinkhuis,
  {\em Normal integral bases and the Spiegelungssatz of Scholz},
  Acta Arith. {\bf 69}  (1995), 1--9
  %

\bibitem{Bush2}  M.R. Bush,
{\em Computation of the Calois groups associated to the $2$-class 
     towers of some quadratic fields}, 
J. Number Theory {\bf 100} (2003), 313--325 
%

\bibitem{Bush} M.R. Bush,
{\em $p$-class towers of imaginary quadratic fields},
thesis Univ. Illionois 2004 
%

\bibitem{BM}  M.R. Bush, D.C. Mayer,
{\em $3$-class field towers of exact length $3$},
J. Number Theory {\bf 147} (2015), 766--777 
%

\bibitem{Dellac} H. Dellac,{\em Question 49},
Interm\'ed. Math. {\bf 1} (1894), 20; ibid. 162--164 
%

\bibitem{Ebb} D. Ebbinghaus (mit V. Peckhaus),
  {\em Ernst Zermelo. An Approach to his Life and Work},
  Springer-Verlag 2007
  %

\bibitem{FLR} G. Frei, F. Lemmermeyer, P. Roquette (Hrsg.), 
  {\em Emil Artin and Helmut Hasse. The correspondence 1923--1958},
  Contributions in Mathematical and Computational Sciences 5;
  Springer 2014
  %
  
\bibitem{FrCE} A. Fröhlich,
  {\em Central Extensions, Galois groups, and ideal class groups of
    number fields}, AMS 1983
  %

\bibitem{FroeUn} A. Fröhlich,
{\em On non-ramified extensions with prescribed Galois group},
Mathematica {\bf 9} (1962), 133--134 
%


\bibitem{FurtD} Ph. Furtwängler,
{\em Zur Theorie der in Linearfaktoren zerlegbaren ganzzahlingen 
     ternären kubischen Formen}, Dissertation Göttingen, 1896 
%

\bibitem{Fw16} Ph. Furtwängler, 
  {\em Über das Verhalten der Ideale des Grundkörpers im
    Klassenkörper}, Monatsh. f. Math. {\bf 27} (1916), 1--15  
%

\bibitem{Fuapp} Ph. Furtwängler,
  {\em Über die simultane Approximation von Irrationalzahlen},
  Math. Ann. {\bf 96} (1927), 169--175
  %

\bibitem{Garb1} D.A. Garbanati,
  {\em The Hasse norm theorem for non-cyclic extensions of the rationals},
  Proc. Lond. Math. Soc. {\bf 37} (1978), 143--164 
  %

\bibitem{Garb2} D.A. Garbanati,
{\em The Hasse norm theorem for $\ell$-extensions of the rationals},
Number Theory and Algebra, 1977, 77--90 
%

\bibitem{Ger3} F. Gerth,
  {\em On $3$-class groups of pure cubic fields},
  J. Reine Angew. Math. {\bf 278/279} (1975), 52--62
  %

\bibitem{Gerth} F. Gerth,
  {\em On $3$-class groups of certain pure cubic fields},
  Bull. Austral. Math. Soc. {\bf 72} (2005), 471--476  
  %

\bibitem{GrLe} K. Grant, J. Leitzel,
  {\em Norm limitation theorem of class field theory},
  J. Reine Angew. Math. {\bf 238} (1969), 105--111
  %

\bibitem{Guy} R.K. Guy,
  {\em Unsolved problems in number theory},
  Springer-Verlag 1981; 3rd ed. 2004
  %

\bibitem{Hajir1} F.~Hajir,
  {\em  On the growth of $p$-class groups in $p$-class field towers},
  J. Algebra {\bf 188} (1997), 256--271 
  %

\bibitem{Hajir2} F.~Hajir,
  {\em On the class numbers of Hilbert class fields},
  Pac. J. Math. {\bf 181} (1997), 177--187 
  %

\bibitem{HasseNm} H. Hasse, 
  {\em Beweis eines Satzes und Widerlegung einer Vermutung über
    das allgemeine Normenrestsymbol},
  Nachr. Ges. Wiss. Göttingen, math.-phys. Kl. 1931, 64--69
  %
  
\bibitem{HS82} F.-P.~Heider, B.~Schmithals,
{\em  Zur Kapitulation der Idealklassen in unverzweigten
        primzyklischen Erweiterungen},
 J. Reine Angew. Math. {\bf 336} (1982),  1--25 
%

\bibitem{HeiderK}  F.-P. Heider, 
{\em Zur Theorie der zahlentheoretischen Knoten},
Diss. Köln 1978 
%

\bibitem{Herbinf} J. Herbrand,
{\em Th\'eorie arithm\'etique des corps de nombres de degr\'e infini},
Math. Ann. {\bf 106} (1932), 473--501; ibid. {\bf 108} (1933), 699--717
%

\bibitem{Hil2a} D. Hilbert,
{\em Über die Theorie der relativ-quadratischen Zahl\-kör\-per},
Jber. DMV {\bf 6} (1899), 88--94 
%

\bibitem{Hil2} D. Hilbert,
  {\em Über die Theorie des relativ-quadratischen Zahl\-kör\-pers},
  Math. Ann.  {\bf 51} (1899), 1--127
  %

\bibitem{Hil3} D. Hilbert,
  {\em Über die Theorie der relativ-Abelschen Zahl\-kör\-per},
  Nachr. Ges. Wiss. Göttingen (1898), 377--399;
  Acta Math. {\bf 26} (1900), 99--132
  %

\bibitem{Jau88} J.~F.~Jaulent,
  {\em  L'\'etat actuel du probl\`eme de la capitulation},
  S\'em. Th\'eor. Nombres Bordeaux, (1987/88), {\bf 17}
  %
 
\bibitem{Jehne} W. Jehne,
{\em On knots in algebraic number theory. In memoriam Arnold Scholz},
J. Reine Angew. Math. {\bf 311/312} (1979), 215--254
%

\bibitem{Kell} O.-H. Keller,
{\em Zur Theorie der ebenen Berührungstransformationen. I},
Math. Ann. {\bf 120} (1949), 650--675 
%

\bibitem{KronDA} L. Kronecker,
{\em Die Periodensysteme von Funktionen reeller Variablen},
Sitz.ber. Preuss. Akad. Wiss. Berlin (1884), 1071--1080; 
Werke III.1, p. 33--46 
%

\bibitem{Krull1} {\sc W. Krull}, 
{\em Idealtheorie in unendlichen algebraischen Zahl\-kör\-pern. I, }
Math. Z. {\bf 29} (1928), 42--54 
%

\bibitem{Krull2} {\sc W. Krull}, 
{\em Idealtheorie in unendlichen algebraischen Zahl\-kör\-pern. II, }
Math. Z. {\bf 31} (1930), 527--557 
%

\bibitem{Kuhnt} Th. Kuhnt,
{\em Generalizations of Golod-Shafarevich and applications},
thesis Univ. Illinois 2000 
%

\bibitem{Langm1} F. Langmayr,
{\em Zur simultanen Diophantischen Approximation. I},
Monatsh. Math. {\bf 86} (1978), 285--300
%

\bibitem{Langm2} F. Langmayr,
{\em Zur simultanen Diophantischen Approximation. II},
Monatsh. Math. {\bf 87} (1979), 133--144
%

\bibitem{Lehmer} E. Lehmer,
{\em Rational reciprocity laws}, Amer. Math. Month. {\bf 85} (1978),
467--472 
%

\bibitem{Lcyc3} F. Lemmermeyer,
{\em Class groups of cyclotomic fields. III},
Acta Arith. 2008 
%

\bibitem{LRH} F. Lemmermeyer, P. Roquette (Hrsg.),
{\em Die mathematischen Tagebücher von Helmut Hasse 1923--1935},
Universitätsverlag Göttingen, 2013 
%


\bibitem{LRS} F. Lemmermeyer, P. Roquette (Hrsg.),
  {\em Der Briefwechsel Hasse - Scholz - Taussky},
  Universitätsverlag G\"ottingen 2016; \\
  \url{https://univerlag.uni-goettingen.de/handle/3/isbn-978-3-86395-253-2}
  %
  
\bibitem{LeoSp} H.W. Leopoldt, 
{\em Zur Struktur der $l$-Klassengruppe galoisscher Zahl\-kör\-per},
J. Reine Angew. Math. {\bf 199} (1958), 165--174
%

\bibitem{Lorenz} F. Lorenz, 
{\em Zur Theorie der Normenreste},
J. Reine Angew. Math. {\bf 334} (1982), 157--170
%

\bibitem{Magnus} W. Magnus, 
{\em Beziehung zwischen Gruppen und Idealen in einem speziellen Ring}, 
Math. Ann. {\bf 11} (1935), 259--280 
%

\bibitem{Maire} C. Maire,
{\em Un raffinement du th\'eor\`eme de Golod-Safarevic}, 
Nagoya Math. J. {\bf 150} (1998), 1--11 
%

\bibitem{MaML} C. Maire, C. McLeman,
{\em On $p^2$-ranks in the class field tower problem},
Ann. Math-. Blaise Pascal {\bf 21} (2014), 57--68 
%


\bibitem{May12a} D.C.~Mayer,
{\em The second $p$-class group of a number field},
Int. J. Number Theory {\bf 8} (2012), 471--505 
%

\bibitem{May12b} D.C.~Mayer,
{\em Transfers of metabelian $p$-groups}, 
Monatsh. Math. {\bf 166} (2012), 467--495 
%

\bibitem{May13} D.C.~Mayer,
{\em The distribution of second $p$-class groups on coclass graphs},
J. Th\'eor. Nombres Bordx. {\bf 25} (2013), 401--456 
%

\bibitem{May14} D.C.~Mayer,
{\em Principalization algorithm via class group structure}, 
J. Th\'eor. Nombres Bordx. {\bf 26} (2014), 415--464 
%

\bibitem{Mayer} J. Mayer,
  {\em Die absolut kleinsten Diskriminanten der biquadratischen
    Zahl\-kör\-per},
  Sitz.ber. Akad. Wiss. Wien IIa, {\bf 138} (1929), 733--724  
%

\bibitem{McL} C.W. McLeman,
{\em A Golod-Shafarevich equality and $p$-tower groups},
thesis Univ. Arizona 2008 
%

\bibitem{MT} E. Menzler-Trott,
{\em Gentzens Problem. Mathematische Logik im nationalsozialistischen
     Deutschland},
Birkhäuser 2001 
%

\bibitem{MinDA} H. Minkowski,
{\em Geometrie der Zahlen}, Leipzig 1896 
%

\bibitem{Mouhib} A. Mouhib,
{\em Sur la tour des $2$-corps de classes de Hilbert des corps 
   quadratiques r\'eels}, 
Ann. Sci. Math. Qu\'ebec {\bf 28} (2004), 179--187 
%

\bibitem{Mori36a} M. Moriya,
  {\em Klassenkörpertheorie im Kleinen für die unendlichen algebraischen
    Zahlkörper},
  Journ. Science Hokkaido Univ. {\bf 5} (1936), 9--66
  %

\bibitem{Mori36b} M. Moriya,
  {\em Klassenkörpertheorie im Großen für unendliche algebraische
    Zahl\-kör\-per},
  Proc. Imp. Acad. Tokyo {\bf 12} (1936), 322--324  
  %

\bibitem{Naito} H. Naito,
  {\em On $\ell^3$-divisibility of class numbers of $\ell$-cyclic extensions},
  Algebraic number theory, Proc. Symp. RIMS, Kyoto/Jap. 1986,
  RIMS Kokyuroku 603, 87--92 (1987)
  %

\bibitem{Neb89} B.~Nebelung,
{\em  Klassifikation metabelscher $3$-Gruppen mit
        Faktorkommutatorgruppe vom Typ $(3,3)$ und Anwendung
        auf das Kapitulationsproblem},
 Diss. Univ. Köln, (1989) 
 %
 
\bibitem{NeissH} F. Neiß,
{\em Darstellung relativ Abelscher Zahl\-kör\-per durch Primkörper 
     und Einheitskörper},
J. Reine Angew. Math. {\bf 166} (1931), 30--53
%

\bibitem{NomGal} A. Nomura, 
{\em On the existence of unramified $p$-extensions with prescribed 
     Galois group},
Osaka J. Math. {\bf 47} (2010), 1159--1165 
%

\bibitem{Nover} H. Nover,
{\em Computation of the Galois groups of $2$-class towers},
thesis Univ. Wisconsin 2009
%

\bibitem{OpGesch} H. Opolka,
{\em Geschlechter von zentralen Erweiterungen},
Arch. Math. {\bf 37} (1981), 418--424 
%


\bibitem{OpoKn} H. Opolka,
{\em Zur Auflösung zahlentheoretischer Knoten},
Math. Z. {\bf 173} (1980), 95--103 
%

\bibitem{OpoKn2} H. Opolka,
{\em Zur Auflösung zahlentheoretischer Knoten in Galois\-er\-wei\-terungen 
     von $\Q$},
Arch. Math. {\bf 34} (1980), 416--420 
%

\bibitem{OpSchur} H. Opolka,
{\em Der Schur-Multiplikator in der algebraischen Zahlentheorie},
Abh. Braunschw. Wiss. Ges. {\bf 33} (1982), 189--195 
%

\bibitem{Opolka} H. Opolka, 
{\em Normenreste in relativ abelschen Zahl\-kör\-per\-er\-wei\-te\-run\-gen
     und symplektische Paarungen},
Abh. Math. Sem. Univ. Hamburg {\bf 54} (1984), 1--4 
%

\bibitem{Ori76} B.~Oriat,
{\em  Spiegelungssatz}, Publ. Math. Fac. Sci. Besan\c{c}on 1975/76
%

\bibitem{Ori79} B.~Oriat,
{\em  Generalisation du 'Spiegelungssatz'},
 Ast\'erisque {\bf 61} (1979),  169--175 
%
 
\bibitem{OS79} B.~Oriat, P.~Satg\'e,
{\em  Un essai de generalisation du 'Spiegelungssatz'},
 J. Reine Angew. Math. {\bf 307/308} (1979),  134--159 
 %
 
\bibitem{Ozaki} M. Ozaki,
{\em Construction of maximal unramified $p$-extensions with prescribed 
     Galois groups},
Invent. Math. {\bf 183} (2011), 649--680 
%

\bibitem{PerronDA} O. Perron,
{\em Über diophantische Approximationen},
Math. Ann. {\bf 83} (192?), 77--84 
%

\bibitem{PoRa} T.P. Pollio, A.S. Rapinchuk,
{\em The multinorm principle for linearly disjoint Galois extensions},
J. Number Theory {\bf 133} (2013), 802--821 
%


\bibitem{ReiKZk} H. Reichardt, 
  {\em Konstruktion von Zahl\-kör\-pern mit gegebener Galoisgruppe von
    Primzahlpotenzordnung},
  J. reine angew. Math. {\bf 177} (1937), 1--5  
  %

\bibitem{RemmFB} V. Remmert,
  {\em Vom Umgang mit der Macht: das Freiburger mathematische
    Institut im ``Dritten Reich''},
  Zeit. f. Sozialgeschichte 20. u. 21. Jhdts. {\bf 14} (1999), 56--85
  %

\bibitem{RemmG} V. Remmert,
  {\em Griff aus dem Elfenbeinturm. Mathematik, Macht und
    Nationalsozialismus: das Beispiel Freiburg},
  Mitt. DMV {\bf 7} (1999), No. 3, 13--24
  %
  
\bibitem{RemmDS} V. Remmert,
  {\em Mathematicians at war. Power struggles in Nazi Germany's mathematical
    community: Gustav Doetsch and Wilhelm S\"uss},
  Revue d'histoire des math\'ematiques 5 (1999), 7--59
  %
  
\bibitem{Rina} H. Richter, 
{\em Über die Lösbarkeit einiger nicht-Abelscher Einbettungsprobleme}, 
Math. Ann. {\bf 112} (1935), 69--84 
%

\bibitem{Riab} H. Richter, 
{\em Über die Lösbarkeit des Einbettungsproblems für 
     Abelsche Zahl\-kör\-per},
Math. Ann. {\bf 112} (1936), 700--726 
%


\bibitem{Sat76} P.~Satg\'e,
{\em  In\'egalit\'es de miroir},
 Sem. Delange-Pisot-Poitou (1967/77), {\bf 18}  4pp 
 %
 
\bibitem{Schoene} Th. Schönemann,
{\em Theorie der symmetrischen Functionen der Wurzeln einer Gleichung. 
     Allgemeine Sätze über Congruenzen nebst einigen Anwendungen 
     derselben}, J. Reine Angew. Math. {\bf 19} (1839), 289--308
%

\bibitem{Sch07} I. Schur,
{\em  Untersuchungen \"{u}ber die Darstellungen der endlichen
       Gruppen durch gebrochene lineare Substitutionen},
 J. Reine Angew. Math. {\bf 132} (1907),  85--137 
 %

\bibitem{Segal} S.L. Segal,
{\em Mathematicians under the Nazis},
Princeton Univ. Press 2003 
%

\bibitem{Ser} J.P. Serre,
{\em Sur une question d'Olga Taussky},
J. Number Theory {\bf 2} (1970), 235--236 
%

\bibitem{Shaf} I.R. Shafarevich,
{\em Construction of fields of algebraic numbers with given 
     solvable Galois group},
Izv. Akad. Nauk. SSSR {\bf 18} (1954), 525--578;
Transl. Amer. Math. Soc. {\bf 4} (1956), 185--237;
Collected Math. Papers 139--237
%

\bibitem{Steinke} G. Steinke, 
{\em Über Auflösungen zahlentheoretischer Knoten},
Schriftenreihe Münster, 1983. 116 pp 
%

\bibitem{Steurer} A. Steurer,
{\em On the Galois groups of the $2$-class field towers of some imaginary
     quadratic number fields},
thesis Maryland 2006 
%

\bibitem{Suzuki} H. Suzuki,
{\em A generalization of Hilbert's theorem 94},
Nagoya Math. J. {\bf 121} (1991), 161--169 
%

\bibitem{Tannaka} T. Tannaka,
{\em Über die Konstruktion der Galoisschen Köprer mit vor\-ge\-ge\-be\-ner
     $p$-Gruppe}, 
T\^ohoku math. J. {\bf 43} (1937), 252--260
%


\bibitem{T52a} O. Taussky, 
{\em Arnold Scholz zum Gedächtnis},
Math. Nachr. 7, 379-386 (1952)
%

\bibitem{Uchi} K. Uchida,
{\em Unramified extensions of quadratic number fields. II},
T\^ohoku Math. J. {\bf 22} (1970), 220--224 
%

\bibitem{WeilA} A. Weil, 
{\em L'avenir des math\'ematiques. Les grands courants de la
      pens\'ee math\'ematique}, 
Marseille 1948 
%

\end{thebibliography}
\end{document}